\documentclass[3p,10pt]{elsarticle}

\journal{}
\date{}
\usepackage{amsmath,amsfonts,amssymb}
\usepackage{graphicx}
\usepackage{hyperref}

\usepackage[svgnames]{xcolor}
\hypersetup{colorlinks=true}
\usepackage[absolute]{textpos}
\TPGrid[0pt,0pt]{20}{30}
\TPMargin{2.0mm}
\TPoptions{showboxes = true }
\textblockcolour{SkyBlue!60}
\textblockrulecolour{NavyBlue}

\graphicspath{{Images/}}

\begin{document}
\hypersetup{ citecolor=teal, filecolor=magenta, urlcolor=NavyBlue }

\begin{frontmatter}

\title{An ALE residual distribution scheme for the unsteady Euler equations over triangular grids with local mesh adaptation}

\author[aff1]{Stefano Colombo\fnref{fn1}}
\ead{stefano.colombo@upm.es}

\author[aff2]{Barbara Re\corref{corrauthor}\fnref{fn1}}
\ead{barbara.re@polimi.it}

\affiliation[aff1]{organization={ETSIA-UPM (School of Aeronautics - Universidad Politécnica de Madrid)},
             addressline={Plaza de Cardenal Cisneros 3},
             city={Madrid},
             postcode={28040},
             country={Spain} }
             
\affiliation[aff2]{organization={Department of Aerospace Science and Technology, Politecnico di Milano},
             addressline={via La Masa 34},
             city={Milano},
             postcode={20156},
             country={Italy} }

\cortext[corrauthor]{Corresponding author}
\fntext[fn1]{This work was initiated while the authors were at Institute of Mathematics, University of Z\"{u}rich}

\begin{abstract}
This work presents a novel interpolation-free mesh adaptation technique 
for the Euler equations within the arbitrary Lagrangian Eulerian framework. For the spatial discretization, we consider a residual distribution scheme, which provides a pretty simple way to achieve high order accuracy on unstructured grids. 
Thanks to a special interpretation of the mesh connectivity changes as a series of fictitious continuous deformations, we can enforce by construction the so-called geometric conservation law, which helps to avoid spurious oscillations while solving the governing equations over dynamic domains.
This strategy preserves the numerical properties of the underlying, fixed-connectivity scheme, such as conservativeness and stability, as it avoids an explicit interpolation of the solution between different grids.
The proposed approach is validated through the two-dimensional simulations of steady and unsteady flow problems over unstructured grids. 
\end{abstract}


\end{frontmatter}

%\linenumbers

%%% For the arXiv version
\begin{textblock}{15}(2.5,1)\large
\textit{Accepted version}

\textbf{Computers \& Fluids (2022) 239:105414} \normalsize
\href{https://doi.org/10.1016/j.compfluid.2022.105414}{DOI: 10.1016/j.compfluid.2022.105414}

The final publication is available at \href{https://doi.org/10.1016/j.compfluid.2022.105414}{sciencedirect.com}
\end{textblock}

\section{Introduction}
In computational fluid dynamics (CFD), mesh adaptation is widely considered an invaluable tool to deal with flow structures that move through the computational domain, with arbitrary large deformations of the boundaries, or with fields characterized by different spatial scales.
In these situations, \textit{h}-adaptation, that is grid connectivity changes such as node insertion, node deletion, or edge swapping, is crucial to comply with the flow evolution and/or with the boundary movement while preserving a good mesh quality, thus to preserve the accuracy of the numerical solution~\cite{Mitchell2014,Alauzet2016,Park2016}.

The goal of mesh adaptation is to optimize a given grid according to an error estimator and a given set of constraints, e.g., minimum or maximum element size, limit on the length ratio between consecutive edges, accuracy of boundary representations, and also computational resources.
A broad palette of mechanisms are available to locally modify the grid: mesh enrichment by means of edge or element split and node insertion by Delaunay triangulation, as well as edge swap and collapse for mesh coarsening, and mesh smoothing algorithm at fixed connectivity ~\cite{Webster1994,Dolejsi1998,Dapogny2014}.
In like manner, several strategies have been proposed to build the error estimator that drives mesh modification, ranging from the simplest ones based on the derivatives of solution variables~\cite{Choi2009,Re2017jet} to more complex ones based on error analyses~\cite{castro1997,Borouchaki1997} or adjoint-based techniques~\cite{Fidkowski2011}.

In addition to an effective error estimator and a wide set of connectivity modifications that can be used to reach the target grid spacing, a further key element of mesh adaptation strategies in CFD is the mapping of the solution from the original to the adapted mesh.
This is especially true in unsteady simulations, when flow features develop and/or the computational domain deforms as time evolves,
or when adaptation is performed in a loop to reach a steady solution: in these situations, we would like to keep track of the solution computed in the previous step and transfer it to the new grid.
As an illustration, consider a given grid $\Omega_h^n$ with $N_j^n$ nodes, over which the solution $\mathbf{U}^n$ has been computed. Now, it undergoes solution-based mesh adaptation which generates a new grid $\Omega_h^{n*}$, characterized by $N_j^{n*}$ nodes, with generally  $N_j^{n*}\neq N_j^n$.
The common strategy to map the solution $\mathbf{U}^n$ over $\Omega_h^{n*}$ is interpolation~\cite{Alauzet2007,Dobrzynski2008,Balan2020}, however this operation may affect negatively the numerical properties of the scheme used for the integration of the partial differential equation (PDEs), e.g., stability and conservativeness, or it may introduce spurious oscillations~\cite{Alauzet2016,Barral2017}.
 
Some attempts to preserve conservation and avoid an explicit interpolation of the solution between different meshes were made by exploiting the arbitrary Lagrangian Eulerian (ALE) formulation, which, in its standard version, aims to solve the governing equations over moving or deforming control volumes, but keeping fixed the grid connectivity.
Mainly, two different paths were proposed to extend the ALE formulation to adaptive meshes: 
a first one is based on a global re-meshing step followed by the construction of new Voronoi tessellation, and a second one focuses on local connectivity changes.
Following the former path, two seminal works were the \textit{Reconnection ALE} (ReALE) algorithm~\cite{Loubere2010}, which kept constant the number of cells between the original and the new mesh to interpolate the solution in a conservative way,
and the \textit{AREPO} code~\cite{Springel2010}, which, during the simulation, adapted the Voronoi tessellation by integrating continuously the motion of a set of mesh-generator points.
Following the second path, Guardone and co-workers~\cite{Guardone2011,Isola2015} proposed a strategy to describe the local grid topology changes due to unstructured mesh adaptation as a series of fictitious collapses and expansions which can be treated within the underlying ALE formulation as continuous deformations, preserving mesh conservation and avoiding any interpolation. This strategy was also extended to three-dimensional grids~\cite{Re2017} and non-ideal fluid flows~\cite{Re2019}.
Finally, still following the local path, also Barral and co-workers~\cite{Barral2017} proposed a three-dimensional adaptive ALE scheme for CFD simulations, which however guarantees conservativeness only if the volume domain does not change.

All the mentioned pioneering works considered, for the spatial discretization of the governing equations, finite volume (FV) schemes that can be at most second-order accurate.
Compared to conventional second-order schemes, high-order methods can increase considerably the solution accuracy for the same computational cost or, equivalently, they can achieve the same accuracy on a coarser mesh, i.e. at a lower cost~\cite{Wang2013}.
Up to our knowledge, the unique extension of an adaptive ALE scheme towards high-order discretization  has been recently proposed by Gaburro et al.~\cite{Gaburro2020}, who have extended the concepts behind \textit{AREPO}  to a discontinuous Galerkin (DG) scheme over two-dimensional unstructured meshes, which are re-generated at each time step.
Although this strategy provides excellent results in 2D, its extension to three-dimensional problems seems to be complex and tightly linked to the adopted space-time integration scheme.
Thereupon, we undertook the study and development of an arbitrary high-order, interpolation-free ALE scheme for unsteady problems over adaptive grids that, following the second of the paths mentioned above, handles local grid modifications.
In this work, we present the first step towards this long-term research goal, which consists in a conservative residual distribution (RD) scheme for two-dimensional adaptive grids able to describe node insertion, node deletion and edge swap within the ALE framework.

In analogy with~\cite{Isola2015,Re2017}, the key role of the proposed strategy is played by a series of fictitious continuous deformations of the mesh elements which allow a straightforward definition of the grid velocities which fulfill the geometric conservation law (GCL) even in presence of topology changes. The GCL is a constraint generally enforced while solving evolutionary equations over moving grids to enhance stability~\cite{Lesoinne1996,Formaggia2004}.
Generally speaking, it requires that the numerical scheme used to solve the governing equations preserves a uniform flow, and it can be fulfilled by a proper computation of the geometric quantities involved in the mesh movement~\cite{Farhat2001}.

Residual distribution (or fluctuation splitting) schemes~\cite{Abgrall2006b} are based on a continuous nodal finite element approximation of the solution and offer the possibility to achieve arbitrarily high order of accuracy on unstructured grids in a pretty simple way~\cite{Abgrall2011,Deconinck2017}.
This type of spatial discretization relies on a compact stencil, an important feature that simplifies the derivation and implementation of the series of fictitious deformations used to describe connectivity changes.
Moreover, the possibility to re-write several standard discretizations under the RD formalism~\cite{Abgrall2006} allows analogies with the interpolation-free adaptive FV schemes that have inspired this work~\cite{Isola2015,Re2017}.
Indeed, we derive the ALE formulation of the Euler equations and the corresponding GCL expression following the aforementioned analogies, although it might be possible to apply the same concept starting from previously proposed ALE RD formulations, as~\cite{Arpaia2015}.
For time integration, we use a backward differentiation formula (BDF), a standard technique that avoids stringent limitations on the time step and easily provides high-order accuracy, to match the same order of accuracy used in the spatial discretization.

The proposed numerical method is implemented in \textsf{Flowmesh}, an inviscid flow solver over dynamic grids for aerodynamic applications under development of the Department of Aerospace Science and Technology, at Politecnico di Milano~\cite{Re2017aiaa,Re2016,IsolaPhd,Cirrottola2021}, which is linked to the open-source library  \textsf{Mmg}~\cite{mmgWeb,Dapogny2014,Dobrzynski2008} for simplicial remeshing.
The validity of the resulting numerical tool is assessed in benchmark tests and in two steady problems, in which mesh adaptation is used to improve the resolution of flow features in the steady state. In these tests, the ALE formulation is simply used to straightforwardly map the solution from one grid to the next one in the loop \textit{mesh adaptation--solution computation}. 
Then, the proposed adaptive ALE RD scheme is applied to simulate the transient evolution of the flow in the forward facing step problem with a fixed geometry, and the unsteady flow around a pitching airfoil. In the latter test, the mesh is both deformed and adapted at every physical time step, to comply with the new position of the solid boundary representing the airfoil and with the time evolution of the flow features.

This paper is divided into further four sections.
In the next one, we describe the development of the proposed GCL-compliant RD scheme for moving grids. We first recall the key ingredients of a standard RD scheme in Sec.~\ref{ss:rdgeneral}, then, in Sec.~\ref{secFV2RD} we define the ALE RD scheme following an analogy with a standard ALE FV scheme. In Sec.~\ref{ss:gcl}, we propose a definition of the interface velocities based on the area swept by the edges during the grid deformation, that fulfils, by construction, the discrete GCL.
In the third section, we explain how we apply the GCL-compliant definition of the interface velocities when a node is inserted or deleted or an edge is swapped by the local mesh adaptation algorithm.
In Sec.~\ref{ss:compProc}, we give an overview of the whole computational procedure.
Sections~\ref{s:validation} and \ref{s:results} present the numerical results of benchmark tests, steady and unsteady problems.
Finally, we discuss the conclusions along with the planned future activities in Sec.~\ref{s:concl}.

\section{An ALE RD scheme for dynamic meshes}\label{s:ale}
The Euler equations, which govern the dynamics of inviscid compressible flows, can be expressed by
\begin{equation}\label{e:eulereqs}
\frac{\partial \textbf{\textit{w}}}{\partial t} +\boldsymbol \nabla \cdot  \mathbf{f}(\textbf{\textit{w}}) = \mathbf{0} ,
\quad \mathrm{for} \;  \mathbf{x}\in \Omega \subset\mathbb{R}^2, t>0 , 
\end{equation}
where the $\textbf{\textit{w}}=[\rho, \, \rho \mathbf{u}, \, E]^\mathrm{T}$ is the vector of conservative variables (density, momentum, and total energy per unit of volume) and 
$\mathbf{f}(\textbf{\textit{w}})$ is the advective flux function, i.e.,
\begin{equation*}
\mathbf{f}(\textbf{\textit{w}})=\begin{bmatrix}
\rho \mathbf{u} \\
\rho \mathbf{u} \otimes \mathbf{u} + P\mathbb{I}_2\\
(E + P) \mathbf{u}
\end{bmatrix}
\end{equation*}
where $P$ is the pressure, and $\mathbb{I}_2$ is the $2 \times 2$ identity matrix.
Appropriate initial and boundary conditions complement the set of PDEs~\eqref{e:eulereqs}, which is closed by the chosen thermodynamic model. Here, we adopt the polytropic ideal gas, so
\begin{equation*}
P=\rho R T \,, \qquad \text{and} \qquad
e=\frac{R T}{\gamma - 1}\,,
\end{equation*}
where $e$ is the internal energy per unit of mass,  $R$ the specific gas constant, $T$ is the temperature, and $\gamma$ is the ratio of specific heats.

In the following subsections, we describe how these governing equations are discretized and solved over a deforming mesh.

\subsection{Steps to define an RD discretization}\label{ss:rdgeneral}
\begin{figure}
    \centering
    \includegraphics[scale=1]{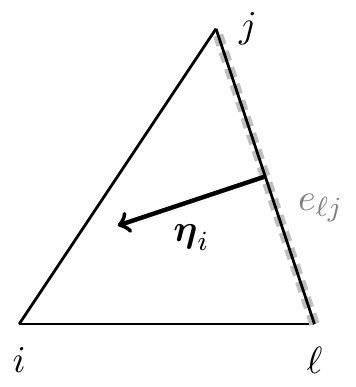}
    \caption{Representation of the normal vector $\boldsymbol{\eta}_i$ associated to the node $i$ of the triangle $K$. The edge opposite to the node $i$ is $e_{\ell j}$, i.e., the edge connecting node $\ell$ to node $j$, highlighted by the gray dashed line.
    The vector $\boldsymbol{\eta}_i$ is defined as the integral of the normal to $e_{\ell j}$, pointing inwards the triangle $K$, over the edge length.
    }
    \label{fig:Eta}
\end{figure}

The domain $\Omega$ is represented by the computational grid $\Omega_h$, composed of the union of $N_K$  non-overlapping triangular elements $K$.
According to the finite element approximation underlying the RD scheme, 
the solution $\textbf{\textit{w}}(\mathbf{x})$ is represented by a globally continuous approximation $\mathit{w}(\mathbf{x})$,
built starting from the nodal values $\mathbf{w}_i=\mathit{w}(\mathbf{x}_i)$.
Hence, given the space of globally continuous piecewise polynomials of degree $k$, i.e. $V_h = \{v \in L^2(\Omega_h) \subset C^0(\Omega_h), v\vert_K \in \mathbb{P}^k, \forall K \in \Omega_h \}$, the solution is written as the linear combination of shape functions $\psi_i \in V_h$:
\begin{equation*}
    \textit{w}(\mathbf{x})=\displaystyle{\sum_{i \in \Omega_h}} \psi_i(\textbf{x}) \mathbf{w}_i = \displaystyle{\sum_{K \in \Omega_h}} \displaystyle{\sum_{i \in K}} \psi_i(\mathbf{x}) \mathbf{w}_i
\end{equation*}
where $i$ are the degrees of freedom (DoFs).
In this work, we consider a linear basis function (so the DoFs coincide with the mesh nodes), which satisfies~\cite{Deconinck2017}
\begin{equation}\label{e:basisfunc}
\psi_j(\mathbf{x}_i)=\delta_{ij} \quad \forall i,j \in \Omega_h;
\quad
\boldsymbol \nabla \psi_i \vert_K = \frac{ \boldsymbol{\eta}_i}{2 |K|}
\end{equation}
where $\delta_{ij}$ is Kronecker delta, $\boldsymbol{\eta}_i$ is the vector normal to the edge of $K$ opposite to the node $i$, pointing towards the node $i$, integrated over the edge length, and $ |K|$ is the area of the triangle $K$ (see Fig.~\ref{fig:Eta}).

%% -----------FIGURE 1----------------------
%% -----------END FIGURE 1-------------------

Now, we briefly recall the main steps required to obtain the spatial RD discretization of~\eqref{e:eulereqs} 
(for the details see, e.g.,~\cite{Abgrall2003,Abgrall2006,Deconinck2017}).
\begin{enumerate}
\item Compute $\forall K \in \Omega_h$ the total residual 
    \begin{equation*}
        \phi^K = \int_K \nabla  \cdot \textbf{f}(\textbf{\textit{w}}) \; \mathrm{d} K  = 
 \int_{\partial K} \!\!\!\! \mathbf{f}(\textbf{\textit{w}}) \cdot\hat{\mathbf{n}} \; \mathrm{d}\partial K,
    \end{equation*}
    which in general is not-null because of the numerical approximation, and where 
$\hat{\mathbf{n}}$ is the outward normal vector of unitary length%
%%% -----Normal unitary-vector
\footnote{Throughout the manuscript, $\hat{\mathbf{n}}$ indicates a unit normal vector, whereas $\boldsymbol{\eta}$ denotes a normal vector integrated over a certain length, defined by the subscript:
for $\boldsymbol{\eta}_i$ the length is the edge of the triangle opposite to $i$, for $\boldsymbol{\eta}_{ij}$ it is the portion of the finite volume interface $\partial \mathcal{C}_{ij}$, as described later. The same subscript-convention is used to distinguish the unit normal vector to these two curves.
In addition, $\hat{\mathbf{n}}$ without subscript is used in flux integrals of a generic variable $Q$, e.g., $\int_\ell Q \cdot \hat{\mathbf{n}} \mathrm{d} \ell$, to denote the outward-facing unit normal vector to the region bounded by $\ell$, and in
boundary conditions to denote the normal pointing inward the domain.}.
%%%%------
In unsteady problems, the total residual includes also the discretized temporal derivative term, as shown in Sec.~\ref{s:time}.
\item Define the nodal residuals $\phi_i^K$ so that conservation property holds 
    \begin{equation}\label{e:consproperty}
        \phi^K = \displaystyle{\sum_{i \in K}} \phi_i^K, \;\;\;\; \forall K \in \Omega_h.
    \end{equation}
    This splitting characterizes the properties of the scheme and it is computed by defining the distribution coefficients $\beta_i$, so that $ \phi_i^K = \beta_i \phi^K$, 
    where due to conservation $\sum_{i \in K}\beta_i = 1$.
 \item The final equation for each DoF is obtained by summing all $\phi^K_i$ for all the elements sharing the node $i$
    \begin{equation*}
        \displaystyle{\sum_{K | i \in K}} \phi_i^K = \mathbf{0} \quad \forall i \in \Omega_h \,,
    \end{equation*}
   which gives rise to a non-linear system of equations, which is solved in a pseudo-time step fashion by looking for the limit solution for $ n \rightarrow \infty$ of
   \begin{equation}\label{e:finalsystemRD}
       \mathrm{C}_i\frac{\mathbf{w}_i^{n+1} - \mathbf{w}_i^{n}}{\Delta \tau} +
       \!\! \sum_{K | i \in K} \phi_i^K = \mathbf{0} \quad \forall i \in \Omega_h\,,
   \end{equation}
   where $\mathrm{C}_i$ is the area of the dual cell associated with the node $i$ and $\tau$ is a fictitious time used for the purpose of solving the system.
\end{enumerate}
In the next subsections, we particularize this general procedure for the GCL-compliant ALE formulation of the governing equations~\eqref{e:eulereqs}.

\subsection{From the ALE FV to the RD discretization}\label{secFV2RD}
The goal of this subsection is to re-write the ALE FV discretization of the Euler equations along with their GCL under the RD formalism.  In this way, we can capitalize on the analogy with the formulation proposed by Guardone and co-workers~\cite{Guardone2011,Isola2015, Re2017}.

Let us consider the first-order FV approximation of the solution over the dual mesh of $\Omega_h$.
With a slight abuse of notation, we denote again as $\mathbf{w}_i$ this piece-wise constant representation of the solution evaluated at node $i$. 
A standard ALE FV discretization of system~\eqref{e:eulereqs}, for the finite volume $\mathcal{C}_i$ defined as the dual cell associated to  node $i$, reads
\begin{equation}\label{FV2RDS}
 \frac{\mathrm{d}(\mathrm{C}_i \mathbf{w}_i)}{\mathrm{d}t} + \!\!
 \sum_{j \in \partial \mathcal{C}_{ij} \mid j \neq i}  \!\!\left( \textbf{\textit{F}}_{ij} \cdot \boldsymbol{\eta}_{ij} - \mathbf{w}_{ij} \nu_{ij} \right) = \mathbf{0} \,,
\end{equation}
where $\mathrm{C}_i$ is the area of the finite volume associated to the node $i$ and the summation is performed over all finite volumes $\mathcal{C}_j$ sharing a portion of their interface with $\mathcal{C}_i$, which is denoted as $\partial \mathcal{C}_{ij}$.
Accordingly, $\textbf{\textit{F}}_{ij}$ is the numerical finite volume evaluation of the flux $\mathbf{f}(\textbf{\textit{w}})$ exchanged across $\partial \mathcal{C}_{ij}$ and
$\boldsymbol{\eta}_{ij}=\int_{\partial \mathcal{C}_{ij}} \hat{\mathbf{n}}_{ij} \mathrm{d}\ell $
is the integrated normal, where $\hat{\mathbf{n}}_{ij} $ points outward $\mathcal{C}_i$.
Between the integrated normal $\boldsymbol{\eta}_{ij}$ (defined in the FV context) and the nodal normal vector ${\boldsymbol{\eta}_i}$ introduced in~\eqref{e:basisfunc} (in the RD context), the following geometrical relation holds within the triangle $K$~\cite{Deconinck2017}:
\begin{equation}\label{e:relNormals}
 \frac{\boldsymbol{\eta}_i}{2} = -\sum_{j \neq i} \boldsymbol{\eta}_{ij}\,. \qquad
    \forall j \in K
\end{equation}
The second contribution in the summation in~\eqref{FV2RDS} is the ALE flux, where
\begin{equation*}
\mathbf{w}_{ij}=\frac{1}{2}(\mathbf{w}_{i}+\mathbf{w}_j) \quad \text{and} \quad
\nu_{ij}=\int_{\partial \mathcal{C}_{ij}} \mathbf{v} \cdot \hat{\mathbf{n}}_{ij} \mathrm{d}\ell \,.
\end{equation*}
The latter, $\nu_{ij}$, is the so-called interface velocity, that is the normal component of the grid velocity $\mathbf{v}$ integrated over the interface portion.
In~\cite{Isola2015,Re2017}, this quantity plays the crucial role to fulfil the GCL.

Equation~\eqref{FV2RDS} is recast within the RD framework (cfr. \eqref{e:finalsystemRD}) by adding two null terms:
\begin{equation}\label{FV2RDS_1}
\begin{aligned}
  & \frac{\mathrm{d}(\mathrm{C}_i \mathbf{w}_i)}{\mathrm{d}t} + {\sum_{K \mid i \in K}}{\sum_{j \neq i}} \big( \textbf{\textit{F}}_{ij} \cdot \boldsymbol{\eta}_{ij} - \mathbf{w}_{ij} \nu_{ij}
     - \textbf{F}_i \cdot \boldsymbol{\eta}_{ij} - \mathbf{w}_i \nu_i \big) \\
= &\frac{\mathrm{d}(\mathrm{C}_i \mathbf{w}_i)}{\mathrm{d}t} + {\sum_{K \mid i \in K}} \phi_i^K = 0
\end{aligned}
\end{equation}
with $\textbf{F}_i = \textbf{f}(\textbf{w}_i)$
and $\nu_i$ the interface velocity associated to the node $i$, yet to be defined.

The summation in~\eqref{FV2RDS_1} occurs on a closed surface, i.e. the finite volume cell, thus the first term we added is always zero since $\sum_{K \mid i \in K} \sum_{j \neq i} \boldsymbol{\eta}_{ij} = 0$.
Regarding the second one, which is not null a priori, we impose the following condition on the
RD interface velocity:
\begin{equation}\label{FV2RDS_2}
    {\sum_{K | i \in K}} \nu_i = 0 \,.
\end{equation}
Thanks to this constraint, we can simply define the total residual  $\phi^K = \sum_{i \in K} \phi_i^K$ as
\begin{equation} \label{FV2RDS_0}
    \phi^K = {\sum_{i \in K}} {\sum_{j \neq i}}\left( -\textbf{F}_i \cdot \boldsymbol{\eta}_{ij} - \mathbf{w}_i  \nu_i \right) \,,
\end{equation}
because, due to conservativity, $ \boldsymbol{\eta}_{ij} = - \boldsymbol{\eta}_{ji}$ and $ \nu_{ji} = - \nu_{ij}$, so many terms of $\phi_i^K$ (see~\eqref{FV2RDS_1}) cancel out while summing over all nodes of the triangle $K$.

To clarify why the expression in~\eqref{FV2RDS_0} is equivalent to the evaluation of the integral of the fluxes in~\eqref{FV2RDS}, we need to better define $\nu_i$, the only unknown term in~\eqref{FV2RDS_0}. To find its expression, we derive the GCL condition for the governing equations~\eqref{e:eulereqs} according to the RD discretization.  So, as standard practice, we impose that Eq.~\eqref{FV2RDS_1} preserves a uniform flow, getting
$$ \frac{\mathrm{d}\mathrm{C}_i}{\mathrm{d}t} + {\sum_{K \mid i \in K}} {\sum_{j \neq i}} \left( -\nu_{ij} - \nu_i \right) = 0 \,,$$
then, we split the time derivative of the area of finite volume $\mathcal{C}_i$ among all the associated triangles, getting
\begin{equation*}  
 {\sum_{K \mid i \in K}}  \Bigg( \frac{\mathrm{d} \,|K_i|}{\mathrm{d}t} - {\sum_{j \neq i}} \nu_{ij} - \nu_i \Bigg) = 
 \sum_{K \mid i \in K} \phi_{i, \mathrm{GCL}}^K = 0 \;,
\end{equation*}
where $$|K_i| = \int_K \psi_i(\mathbf{x}) \; dK =\frac{|K|}{3}$$
is the portion of the area of the triangle $K$ pertaining to the node $i$.
Following similar arguments as in the derivation of~\eqref{FV2RDS_0},
the total residual for the GCL $\phi_{\mathrm{GCL}}^K = \sum_{i \in K} \phi_{i, \mathrm{GCL}}^K$ is
\begin{equation*}
    \phi_{\mathrm{GCL}}^K = {\sum_{i \in K}} \left( \frac{\mathrm{d} \,|K_i|}{\mathrm{d}t} - \nu_i \right) \,.
\end{equation*}
A possible definition of the interface velocities that fulfils Eq.~\eqref{FV2RDS_2} is obtained by imposing that all the total residuals just defined are zero, thus $\phi_{\mathrm{GCL}}^K = 0, \quad \forall K \in \Omega_h$ and
\begin{equation}\label{FV2RDS_4}
    {\sum_{i \in K}} \nu_i = \frac{\mathrm{d} \,|K|}{\mathrm{d}t}, \qquad
    \forall K \in \Omega_h\,.
\end{equation}
In summary, to fulfil the GCL in the considered framework, the definition of $ \nu_i$ must satisfy the two conditions stated in~\eqref{FV2RDS_2} and \eqref{FV2RDS_4}.

Now that we have formalized the GCL constraint, we can come back to the total residual for the Euler equations defined in~\eqref{FV2RDS_0}.
First of all, knowing that $\frac{\mathrm{d}\,|K|}{\mathrm{d}t} = \oint_{\partial K} \mathbf{v} \cdot \hat{\mathbf{n}} \; \mathrm{d} \partial K$ and expressing the grid velocity over the element $K$ as 
$\psi_i(\mathbf{x}) \mathbf{v} $, we recast Eq.~\eqref{FV2RDS_4} as
\begin{equation*}
    {\sum_{i \in K}} \nu_i = {\sum_{i \in K}} \oint_{\partial K} \psi_i(\mathbf{x}) \mathbf{v} \cdot \hat{\mathbf{n}} \;\mathrm{d}\partial K.
\end{equation*}
Then, using the geometrical relations given in~\eqref{e:basisfunc} and \eqref{e:relNormals}, and assuming a piece-wise linear interpolation of the flux as
 $\textbf{f} = \sum_{i \in K} \psi_i(\mathbf{x}) \mathbf{F}_i$, we manipulate Eq.~\eqref{FV2RDS_0} as
\begin{align*}
\phi^K &= {\sum_{i \in K}} \Big( -\mathbf{F}_i {\sum_{j \neq i}} \boldsymbol{\eta}_{ij} \Big) - {\sum_{i \in K}} \mathbf{w}_i \nu_i = \\
    &= {\sum_{i \in K}} \; \mathbf{F}_i \frac{\boldsymbol{\eta}_i}{2} - {\sum_{i \in K}} \mathbf{w}_i \oint_{\partial K} \psi_i (\mathbf{x}) \mathbf{v} \cdot \hat{\mathbf{n}} \; \mathrm{d} \partial K \\ 
    &= {\sum_{i \in K}} \Big( \! \int_K \!\! \mathbf{F}_i \nabla \psi_i(\mathbf{x}) \; \mathrm{d}K - \oint_{\partial K} \!\! \mathbf{w}_i  \psi_i (\mathbf{x}) \mathbf{v} \cdot \hat{\mathbf{n}} \; \mathrm{d} \partial K \Big)  \\
    &=\int_K \nabla \cdot \mathbf{f} \; \mathrm{d}K - \oint_{\partial K} \textbf{\textit{w}} \mathbf{v} \cdot \hat{\mathbf{n}} \; \mathrm{d} \partial K \;
    = \int_K \nabla \cdot ( \mathbf{f} - \mathbf{\textit{w}} \mathbf{v} ) \; \mathrm{d} K
\end{align*} 
which proves that eq. \eqref{FV2RDS_0} is equivalent to compute an integral over the element $K$, thus a residual. 

Finally, we need to choose the first-order scheme underlying the residual definition. 
In this work, we use the Lax-Friedrichs scheme, which has a low computational cost. So, following~\cite{Abgrall2006b}, the nodal residuals are defined as
\begin{equation}\label{eq:laxFriedrichs}
    \phi_i^K = \frac{\phi^K}{3} + \alpha_{\mathrm{LF}} \left( \mathbf{w}_i - \overline{\mathbf{w}}^K \right) \,,
    \qquad \text{with} \quad
    \overline{\mathbf{w}}^K=\frac{1}{3}\sum_{j\in K} \mathbf{w}_j \,.
\end{equation}
This splitting satisfies the conservation property~\eqref{e:consproperty}, while the positivity of the scheme is determined by the parameter $\alpha_{\mathrm{LF}}$.
To define it, we consider a generic linearization of the flux
$\mathbf{F}_i \approx  \mathbb{A}(\tilde{\mathbf{w}}) \mathbf{w}_i $, so that the  local residual can be written as
\begin{equation*}
    \phi_i^K = \frac{1}{3} \sum_{j \in K}
    \left( \mathbb{A}(\tilde{\mathbf{w}}) \cdot \frac{\boldsymbol{\eta}_j}{2} - \nu_j \mathbb{I}_4 \right) \mathbf{w}_j + \alpha_{\mathrm{LF}} \left( \mathbf{w}_i - \overline{\mathbf{w}}^K \right) 
      = \frac{1}{3} {\sum_{j \in K}} \left( \mathbb{K}_j - \nu_j \mathbb{I}_4 \right) \mathbf{w}_j + \alpha_{\mathrm{LF}} \left( \mathbf{w}_i - \overline{\mathbf{w}}^K \right)
\end{equation*}
where $\mathbb{K}_j = \mathbb{A}(\tilde{\mathbf{w}}) \cdot \frac{\boldsymbol{\eta}_j}{2}$ and $\mathbb{I}_4$ is the $4\times 4$ identity matrix.
To ensure local positivity, the parameter $\alpha_{\mathrm{LF}}$ should be chosen larger than the spectral radii of the matrices $\mathbb{K}_j - \nu_j \mathbb{I}_4$ for all $j\in K$.
%\begin{equation*}
%    \alpha_{\mathrm{LF}} \ge \rho(\mathbb{K}_j - \nu_j \mathbb{I})
%\end{equation*}
%where $\rho(\mathbb{T})$ is the spectral radius of the matrix $\mathbb{T}$. 
Thus, we define $\alpha_{\mathrm{LF}}$ as 
\begin{equation}\label{alpha}
    \alpha_{\mathrm{LF}} = \max_{j \in |K|} \biggl \lvert  \frac{1}{2} \left( \mathbf{u}(\tilde{\mathbf{w}}) \cdot \boldsymbol{\eta}_j \pm c(\tilde{\mathbf{w}}) \; \lvert \boldsymbol{\eta}_j \rvert  \right) - \nu_j  \biggl \rvert \,,
\end{equation}
where $c$ is the speed of sound, computed according to the polytropic ideal gas as $c(\mathbf{w}_j)=\sqrt{\gamma R T(\mathbf{w}_j)}$.

\subsection{Boundary conditions}
The contribution from the domain boundary $\partial\Omega$ is taken into account through a corrective approach, as described in~\cite{Agrall2014,AbgrallDeSantis2015}.
For a node $i$ lying on the boundary, the total residual is first computed, as described before, without considering the boundary contribution, then it is corrected by adding the boundary residual defined as
\begin{equation}\label{eq:bintegral}
    \phi_i^{K,\partial}
    = \int_{\Gamma^K} \!\!\!\! \psi_i \left( \mathbf{f}(\boldsymbol{w}^{\partial}) - \mathbf{f}(\textbf{w}) \right) \cdot \hat{\mathbf{n}} \, \mathrm{d} \, \partial \Omega   
    = \int_{\Gamma^K} \!\! \psi_i \boldsymbol\delta\mathbf{f} \, \mathrm{d} \, \partial \Omega 
\end{equation}
where
$\boldsymbol {w}^{\partial}$ is the physical state that takes into account the boundary conditions (in a weak sense),
$\mathbf{w}$ is the numerical state, and
$\Gamma^K$ is the union of the edges of the triangle $K$ pertaining to boundary $\partial \Omega$.
The correction flux $\mathbf{f}(\boldsymbol{w}^{\partial}) - \mathbf{f}(\textbf{w})$ vanishes when the boundary conditions impose correctly that the boundary solution equals the state $\boldsymbol{w}^{\partial}$. Now, we define this flux for the boundary conditions used in this work, that is non-penetration and inflow/outflow conditions. 

For the non-penetration condition, considering a solid boundary moving with velocity $\mathbf{v}$, we have:
\begin{equation}
    \delta \mathbf{f} = -v_n \begin{pmatrix} 
    \rho \\
    \rho \mathbf{u} \\
    (E + P)
    \end{pmatrix}
\end{equation}
with $v_n = (\mathbf{u} - \mathbf{v}) \cdot \hat{\mathbf{n}}$ the normal component of the relative velocity between the solid boundary and the flow.

Inﬂow/outﬂow boundary conditions are imposed through a characteristics-based procedure, which automatically determines the number of conditions to be imposed~\cite{Agrall2014}. Considering only the advective flux, it can be linearized at the boundary state as
$$
\mathbf{f}(\boldsymbol{w}^{\partial}) \cdot \hat{\mathbf{n}} \approx
\mathbb{A}^+_\mathbf{n}(\mathbf{w}) \mathbf{w} + \mathbb{A}^-_\mathbf{n}(\mathbf{w}) \boldsymbol{w}^{\partial}
,\quad \text{with} \quad
\mathbb{A}^\pm_\mathbf{n} = \mathbb{R}_\mathbf{n} \Lambda^\pm_\mathbf{n} \mathbb{L}_\mathbf{n}
$$
where $\mathbb{R}_\mathbf{n}$ and $\mathbb{L}_\mathbf{n}$ are, respectively, the matrices of the right and left eigenvectors projected along the direction of $\hat{\mathbf{n}}$,  $\Lambda_\mathbf{n}$ is the diagonal matrix of the eigenvalues,
and the operator $^+$ (or $^-$) selects only the positive (or negative) components, setting to zero the other ones.
Hence, considering also the ALE contribution, the correction flux for an inflow/outflow condition is 
\begin{equation}
    \delta \mathbf{f}= \Big(
    \mathbb{A}_{\mathbf{n}}(\mathbf{w}) - (\mathbf{v} \cdot \hat{\mathbf{n}}) \mathbb{I}_{4}\Big)^{\!-} \left( \boldsymbol{w}^{\partial} - \mathbf{w} \right)\,.
\end{equation}

The boundary residual defined in~\eqref{eq:bintegral} is computed through a quadrature formula.
For the boundary conditions used in this work,
the residual to be added to the node $i$ belonging to the boundary edge $\Gamma^K$, of length $\gamma^K$, becomes, for the non-penetration condition
\begin{equation}
    \phi_i^{K,\partial} = \frac{1}{2}\gamma^K (\mathbf{u}\cdot \hat{\mathbf{n}})_i \begin{pmatrix} \rho \\
    \rho \textbf{u} \\ 
    (E + P) \end{pmatrix}_{\!\!i}  - \nu_i \textbf{w}_i 
\end{equation}
and for the inflow/outflow condition
\begin{equation}
     \phi_i^{K,\partial} = \left( \frac{1}{2} \gamma^K  \mathbb{A}_{\mathbf{n}}(\mathbf{w}_i) - \nu_i \mathbb{I}_{4}\right)^{\!-}
    \left( \boldsymbol{w}^{\partial} - \mathbf{w}_i \right)
    \,.
\end{equation}

\subsection{Time integration}\label{s:time}
As anticipated in Sec.~\ref{ss:rdgeneral}, the temporal term is discretized and added to the total residual $\phi_i^K$. Hence, the resulting space-time residual is
\begin{equation}\label{e:timeint}
   \Phi_i^K = \frac{{\sum_{q = -1}^p}a_q |K_i|^{n-q} \textbf{w}_i^{n-q}}{\Delta t} + \phi_i^K \,,
\end{equation}
where we have used as time integrator a BDF with $p+1$ steps.  
As shown in \cite{Caraeni2006}, this choice provides a straightforward way of getting a second order accuracy in time, when $p=1$.
This space-time residual is inserted in the Lax-Friedrichs scheme defined in~\eqref{eq:laxFriedrichs}, using the same definition of $\alpha_{\mathrm{LF}}$ as in~\eqref{alpha}.

In addition, let us recall that the definition of GCL-compliant interface velocities requires us to take into account the variation of the elements area occurring during the time step $\Delta t$ due to grid adaptation. So we apply the same time integration scheme also for the discretization of the temporal derivative of the element areas in~\eqref{FV2RDS_4}, as explained in the next section, Sec.~\ref{ss:gcl}.
Finally, the role of the time integration in the simulations is better illustrated by the summary of the whole computational procedure in Sec.~\ref{ss:compProc}.

\subsubsection{A note about steady simulations}\label{ss:timesteady}
For validation purposes, we carry out steady simulations, in which we perform a loop \textit{mesh adaptation--solution computation} to optimize the grid according to the steady flow field, without considering any time evolution or transient state.
In such a simulation, we can use the proposed scheme to map the solution from one grid to the next (adapted) one, without any explicit interpolation of the solution.
Admittedly, this is not a mandatory step as the steady computation does not need to keep into consideration the solutions over previous grids, but this clearly improves the convergence, as we can restart the solution computation from an already developed field, rather than the initial, generally uniform, field.
Even in such a steady framework, the GCL-condition requires us to take into account the variation in the element area occurring during the adaptation phase, which is here considered as a \textit{fictitious} time step $\Delta t$.
In this case, the node residual becomes
\begin{equation}
    \Phi_i^K = \frac{|K_i|^{n+1}\textbf{w}_i^{n+1} - |K_i|^n\textbf{w}_i^n}{\Delta t} + \phi_i^K
\end{equation}
where, differently from~\eqref{e:timeint}, we use a simple implicit Euler scheme as the temporal accuracy in this non-physical time does not affect the solution accuracy.

\subsection{Calculation of GCL-complaint interface velocities}\label{ss:gcl}
Now, an open point concerns the computation of the interface velocities $\nu_i$ when the mesh is deformed from $\Omega_h^{n}$ to $\Omega_h^{n+1}$.
It is important to realize that the positions of the grid nodes $\mathbf{x}^{n+1}$ are known, because the solution $\mathbf{U}^{n+1}$ is computed after the grid deformation has been accomplished.
Hence, the problem is to define suitable interface velocities $\nu_i$ that, starting from the available grid information, satisfy the relations derived in Sec.~\ref{secFV2RD} at a discrete level, so that the GCL constraint is fulfilled.
To this end, we recast the BDF scheme in an incremental fashion as $\frac{dy}{dt} \approx \frac{1}{\Delta t} {\sum_{q = -1}^{p-1}} \alpha_q \Delta y^{n-q}$ and use this approximation in Eq.~\eqref{FV2RDS_4} to compute the variation of the area of the triangle $K$  during the time step $\Delta t$. For instance, considering a BDF with 2 steps ($p=1$), Eq.~\eqref{FV2RDS_4} for a triangle $K$ can be integrated in time as
\begin{equation*}
    \frac{1}{\Delta t} \Big( \alpha_{-1} \Delta A^{n+1} + \alpha_{0} \Delta A^{n} \Big)
    =  \sum_{i \in K} \nu_i^{n+1},
    \qquad \text{with}\quad
    \alpha_{-1}=\frac{3}{2} \quad \text{and} \quad \alpha_{0}=-\frac{1}{2}, 
\end{equation*}
where $\Delta A^{n+1} = |K|^{n+1} - |K|^n$ is the variation of the triangle area occurring during the current time step (due to the grid deformation from $\Omega_h^n$ to $\Omega_h^{n+1}$), while
$\Delta A^{n} = |K|^{n} - |K|^{n-1}$ is the variation occurred in the previous time step.
In general, the integration of Eq.~\eqref{FV2RDS_4} with a BDF with $p+1$ steps involves the variations of the triangle area occurred from $p$ time steps ago to the current one:
\begin{equation}\label{e:dgcl}
  \frac{1}{\Delta t} {\sum_{q = -1}^{p-1}} \alpha_q \Delta A^{n-q}
  = \sum_{i \in K} \nu_i^{n+1}  \,, \quad \forall K \in \Omega_h\,.
\end{equation}
Hence, the contribution $\Delta A^{n+1}$ computed in the current time step needs to be stored for the next $p$ steps.

\begin{figure}
    \centering
    \includegraphics[scale=1]{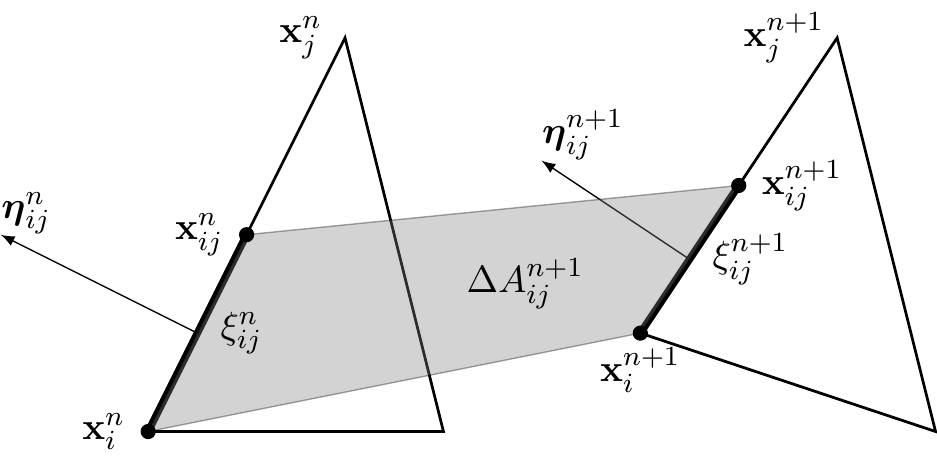}
    \caption{Area swept $\Delta A_{ij}^{n+1}$ by the semi-edge $\xi_{ij}$ during the time step from $t^n$ to $t^{n+1}$. The triangle on the left illustrates the initial position, on $\Omega_h^n$; the one on the right represents the position after the deformation, so on $\Omega_h^{n+1}$;
    $\mathbf{x}_{i}$ and $\mathbf{x}_{j}$ indicate the position of the nodes $i$ and $j$, connected by the edge $e_{ij}$.
    The position of the mid-point of this edge is
    $\mathbf{x}_{ij}= \frac{1}{2} ( \mathbf{x}_{i} + \mathbf{x}_{j})$.
    The arrows illustrate the integrated normal vectors $\boldsymbol{\eta}_{ij}^n$ and $\boldsymbol{\eta}_{ij}^{n+1}$, which differ both in magnitude and direction. The area $\Delta A_{ij}^{n+1}$ amounts to the gray quadrilateral between  $\xi_{ij}^n$ and $\xi_{ij}^{n+1}$ and is computed according to~\eqref{e:sweptarea}.
    }
    \label{f:areaswept}
\end{figure}

Expression~\eqref{e:dgcl}, highlighting the link between the interface velocities $\nu_i^{n+1}$ and the change of the area due to the mesh movement, provides guidance for the computation of GCL-compliant interface velocities.
Now, there could be different ways to fulfil simultaneously Eqs.~\eqref{e:dgcl} and~\eqref{FV2RDS_2}. A simple possibility consists in splitting the change of area of $K$ during the time step among different contributions pertaining to its nodes.
To explain this procedure, which is applied to each element of the grid, let us consider a single triangle $K$ that undergoes a deformation as displayed in Fig.~\ref{f:areaswept}.
Assume to be able to compute the area swept by each semi-edge  $\xi_{ij}= \mathbf{x}_i-\mathbf{x}_{ij}$, where $\mathbf{x}_{ij}$ is the midpoint of the edge $e_{ij}$ between nodes $i$ and $j$. 
This contribution is labelled as $\Delta A_{ij}$ and it is illustrated in Fig.~\ref{f:areaswept}.
Then, we associate to each node $i$ the area $\Delta A_{i} = \sum_{j\in E_i},  \Delta A_{ij}$, where $E_i$ is the set of edges of $K$ connected to $i$.
It is easy to verify that $|K|^{n+1} - |K|^n = \Delta A^{n+1} = \sum_{i \in K} \Delta A_i^{n+1}$.
According to this splitting, we define the GCL-compliant interface velocities on each triangle of the deformed grid $\Omega_h^{n+1}$ as
\begin{equation}\label{e:defnu_i}
    \nu_i^{n+1} = \frac{1}{\Delta t} {\sum_{q = -1}^{p-1}} \alpha_q \Delta A_i^{n-q} \,,
    \qquad \forall i\in K,
\end{equation}
where the area  $\Delta A_{i}^{n-q}$ for $q=0, \dots, p-1$ are known from the previous steps.

If we assume that the grid velocities are constant during the deformation from $\Omega^{n}$ to $\Omega^{n+1}$, the area swept by the semi-edge $\mathbf{x}_i-\mathbf{x}_{ij}$ can be computed as
\begin{equation}\label{e:sweptarea}
  \Delta A_{ij}^{n+1}\!=\! \frac{\left( \mathbf{x}_i^{n+1} \!\! -  \mathbf{x}_{ij}^{n+1}\right) -
  \left( \mathbf{x}_i^{n} -  \mathbf{x}_{ij}^{n}\right)}{4}\cdot
  \left( \boldsymbol{\eta}_{ij}^{n+1} \!\! + \boldsymbol{\eta}_{ij}^n \right)
\end{equation}
where $\boldsymbol{\eta}_{ij} = \int_{\xi_{ij}} \hat{\mathbf{n}} \mathrm{d}\ell = (\mathbf{x}_i - \mathbf{x}_{ij}) \times \hat{\mathbf{n}} $ and $\hat{\mathbf{n}}$ is the outward unit-length vector normal to the edge $e_{ij}$. A graphical definition of all quantities involved in~\eqref{e:sweptarea} is given in Fig.~\ref{f:areaswept}.

Now, let us consider what happens outside the triangle $K$.
Equation~\eqref{e:sweptarea} shows that the area swept by the semi-edge $\xi_{ij}$ of the triangle $K$ has the same magnitude but opposite sign with respect to the one swept by the adjacent triangle $K'$ sharing the nodes $i$ and $j$ with $K$. Indeed, the relation $\boldsymbol{\eta}_{ij}^K = - \boldsymbol{\eta}_{ij}^{K'} $ holds. Thus, $\sum_{K \mid i \in K} \Delta A_i = 0$, and Eq. \eqref{FV2RDS_2} is satisfied.
If the element $K$ lies on the boundary, the contribution of the semi-edge belonging to the boundary is balanced by the area swept by the boundary edge, computed in similar way as in~\eqref{e:sweptarea}, but considering a normal vector pointing toward the inner of the triangle.
Thus, also in this case, the above condition on the integrated normal vectors holds.    

%% -----------FIGURE 2----------------------
%% -----------END FIGURE 2-------------------

\subsection{Building a second-order accurate RD scheme}
The last step of the construction of the proposed ALE RD discretization for dynamic meshes is the extension  toward a  second-order scheme, which is carried out following the procedure described in~\cite{Agrall2014}, and here summarized.
The idea is to start from the first-order scheme for which the coefficients $\beta_i$ are, in general, unbounded, and to apply to them a suitable mapping $\beta_i \rightarrow \hat{\beta}_i$.
Through this mapping, we can obtain bounded coefficients, satisfying conservation and monotonicity constraints, and then compute the high order residuals as $\hat{\Phi}^K = \hat{\beta}_i \Phi_i^K$, taking into account that, as a consequence of the Godunov theorem, a monotonicity and linearity preserving scheme must be non-linear. 

The non-linear mapping is built in the characteristic space, so the first order residual are re-written as
\begin{equation*}
    \Phi_i^{\star,K} = \mathbb{L}_{\textbf{n}} \Phi_i^K \;\;\;\; \mathrm{and} \;\;\;\; \Phi^{\star,K} = {\sum_{i \in K}} \Phi_i^{\star, K}
\end{equation*}
where $\mathbb{L}_{\textbf{n}}$ is the matrix of the eigenvectors projected along the direction of $\textbf{n}$ which can be arbitrarily chosen. The map for high-order is built as follows:
\begin{equation*}
    \hat{\beta}_i = \frac{\max \left( \frac{\Phi_i^{\star, K}}{\Phi^{\star, K}}, 0\right)}{{\sum_{j \in K}} \max_{j \in K} \left( \frac{\Phi_j^{\star, K}}{\Phi^{\star, K}}, 0\right)} \,,
\end{equation*}
and the limited residuals are computed as $\hat{\Phi}_i^{\star, K} = \hat{\beta}_i \Phi^{\star, K}$. The high order residual are projected back into the physical space, getting $\hat{\Phi}_i^K = \mathbb{R}_{\mathbf{n}} \hat{\Phi}_i^{\star, K}$, with $\mathbb{R}_{\mathbf{n}} = \mathbb{L}_{\textbf{n}}^{-1}$.

The high order scheme associated with the Lax-Friedrics first order residual shows spurious oscillation in smooth regions due to the absence of an upwind mechanism. Thus a stabilization term $\Xi$ is required. Following \cite{Ricchiuto2009}, we define for each node $i$ of $K$
\begin{equation*}
    \Xi_i = \frac{\left( \mathbb{A} \cdot \boldsymbol{\eta}_i \right)
    }{2 |K|}  \mathbb{\tau} \Phi^{K} + \frac{1}{12} \displaystyle{\sum_{j \in K}} d_{ij}  (|K_i|^{n+1} \textbf{w}_i^{n+1}-
    |K_i|^{n} \textbf{w}_i^{n})
\end{equation*}
with  $d_{ij}$ the $i,j$ component of $ \mathbb{D} = \footnotesize
\begin{bmatrix}
2 & -1 & -1 \\
-1 & 2 & -1 \\
-1 & -1 & 2 \\
\end{bmatrix}$,
$\mathbb{A}$~the Jacobian matrix of the fluxes, and
$$ \mathbb{\tau} = \dfrac{1}{2} |K| \Big( \displaystyle{\sum_{j \in K}} \widetilde{\mathbb{K}}_j^+ \Big)^{-1}, \quad \text{with} \quad
\widetilde{\mathbb{K}}_j^+ = \mathbb{R}_{\boldsymbol{\eta}_j} \boldsymbol{\Lambda}_{\boldsymbol{\eta}_j}^+ \mathbb{L}_{\boldsymbol{\eta}_j}\,, $$
where $\boldsymbol{\Lambda}_{\boldsymbol{\eta}_j}^+$ is the diagonal matrix of the positive eigenvalues.

Finally, the high order residual can be written as
\begin{equation*}
    \Phi_i^{\mathrm{HO}, K} = \hat{\Phi}_i^K + \theta^K \Xi \,,
\end{equation*}
where $\theta^K$ is a filtering term that acts only in smooth regions of the flow, where $\theta^K \sim 1$, whereas $\theta^K \sim 0$ close to discontinuities. In this work, $\theta^K$ is defined as
\begin{equation*}
%    \theta^{K} = 1 - \max_{j \in K} \left( \max_{K' \mid j \in K'} \max_{j \in K'} \frac{\lvert s_j - \bar{s}_{K'} \rvert}{s_j + \bar{s}_{K'} + \epsilon} \right)
    \theta^K = \min \left( 1, \frac{1}{\frac{\rho_t}{ \lvert K \rvert} + \epsilon} \right)
\end{equation*}
with $\rho_t = v_t \Phi^{K}$ where $v_t$ is an average on the element of the entropy variables and $\epsilon = 10^{-10}$.

\section{Extension of the ALE RD scheme to adaptive grids}\label{s:adapt}
The key concept that paves the wave to the description of mesh connectivity changes within the ALE framework is the possibility to define GCL-compliant interface velocities also when the mesh is locally adapted.
In these situations, we can describe the topological modifications as a series of continuous deformations that occur within the fictitious time step $0\leq \tau \leq 1$.
Thanks to this interpretation, the area swept by the edges of the involved elements can be still computed according to~\eqref{e:sweptarea}, considering as initial one the position of the nodes on the original grid, i.e., $\mathbf{x}^{n}=\mathbf{x}(\tau=0)$, and as final one the position of the nodes on the adapted grid, i.e., $\mathbf{x}^{n+1}=\mathbf{x}(\tau=1)$.

\begin{figure}
\centering
\includegraphics[scale=1]{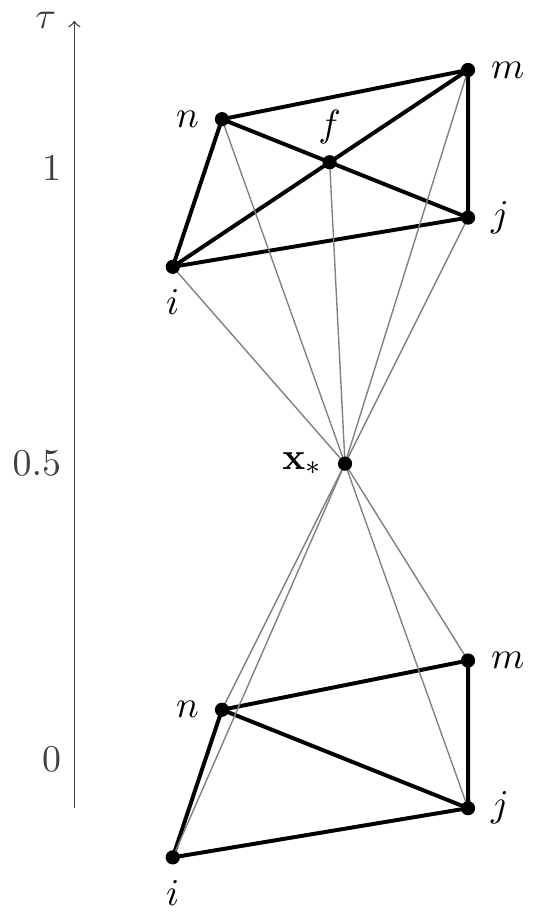}%
    \caption{Collapse-expansion procedure for a point insertion by edge split.
    The evolution of topology change is illustrated  within the fictitious time $0\leq \tau \leq 1$. The initial configuration ($\tau=0$) is shown in the bottom part, moving toward the top $\tau$ increases, until reaching the final configuration, shown in the top part, for $\tau=1$. The gray lines display the movement of the nodes of the triangles involved in the modifications (the set  $\mathcal{K}^\mathrm{ad}$).
    For $\tau=0.5$, the lines collides in one point, $\mathbf{x}_*$ and the elements have null area; in this moment, the connectivity is modified. In this figure, at  $\tau=0.5$, a new node is inserted in the midpoint of the edge $e_{nj}$ shared by the two original triangles; for $\tau>0.5$, a new line shows the position of the added node within $\mathcal{K}^\mathrm{ad}$; in the final configuration, two new elements are present.}
    \label{fig:PointAddition}
\end{figure}

\begin{figure}
\centering
\includegraphics[scale=1]{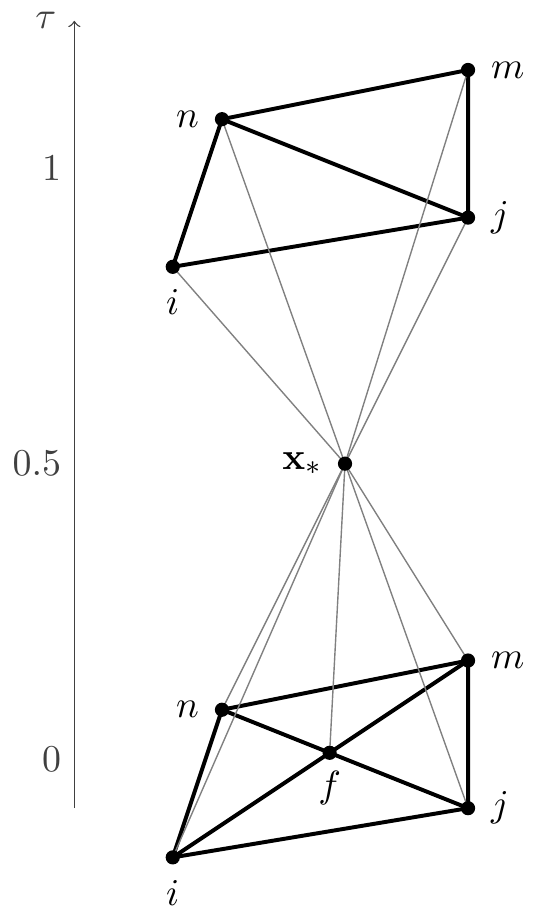}
    \caption{Collapse-expansion procedure for a node deletion.
    Following the same nomenclature used in Fig.~\ref{fig:PointAddition}, the evolution of the triangles involved in the deletion of node $f$ (the set  $\mathcal{K}^\mathrm{ad}$) is illustrated.
    At $\tau=0.5$, the node shared by four triangles initially in $\mathcal{K}^\mathrm{ad}$ is deleted, along with two triangles; in the new configuration for $\tau>0.5$, only two triangles form $\mathcal{K}^\mathrm{ad}$.}
    \label{fig:PointDeletion}
\end{figure}

\begin{figure}
\centering
\includegraphics[scale=1]{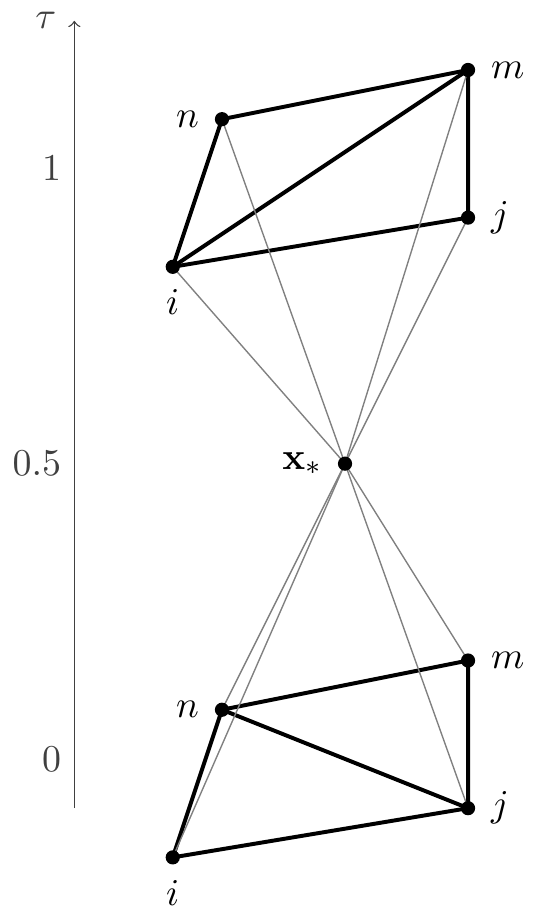}
    \caption{
    Collapse-expansion procedure for an edge swapping.
    Following the same nomenclature used in Fig.~\ref{fig:PointAddition}, the evolution of the triangles involved in the swapping of the edge $e_{nj}$ with the edge $e_{im}$ (the set  $\mathcal{K}^\mathrm{ad}$) is illustrated.
    At  $\tau=0.5$, the edge shared by the two original triangles is swapped to connect the opposite vertices; still two triangles form $\mathcal{K}^\mathrm{ad}$ for $\tau>0.5$, but they have a different topology from the initial ones. For all cases, we highlight that the position of the nodes on the external boundary of the set $\mathcal{K}^\mathrm{ad}$ is exactly the same at  $\tau=0$ and at $\tau=1$. This allows the exclusion from the collapsed-expansion procedure of the edges external to this region, which, even if included, would sweep the same area during the collapse and expansion phase, but with opposite sign, so $\Delta A_{ij} =0 \; \forall e_{ij} \notin \mathcal{K}^\mathrm{ad}$.}
    \label{fig:Edgeswapping}
\end{figure}

This collapse-expansion procedure is illustrated for node insertion, node deletion, and edge swap in Fig.~\ref{fig:PointAddition},~\ref{fig:PointDeletion} and~\ref{fig:Edgeswapping} respectively. In general, when the mesh adaptation algorithm determines to perform a certain topology change, the following steps are followed:
\begin{itemize}
    \item at $\tau=0$: identification of the set of elements involved in the topology change, labeled $\mathcal{K}^\mathrm{ad}$: this set is defined so that the edges of all elements $K^\prime \notin \mathcal{K}^\mathrm{ad}$ are not affected by topological modifications;
    \item $0 < \tau < 0.5$: all triangles $K \in \mathcal{K}^\mathrm{ad}$ are collapsed over a point $\mathbf{x}_*$, which lies within $\mathcal{K}^\mathrm{ad}$, and during this deformation we evaluate the area swept by all semi-edges and compute the nodal contributions $\Delta A_i^\mathrm{col}\quad \forall i \in \mathcal{K}^\mathrm{ad}$ (where the superscript ``col'' stands for \textit{collapse});
    \item at $\tau=0.5$: all elements in $\mathcal{K}^\mathrm{ad}$ have null area and zero-length edges, so we actually modify the topology without affecting the solution as all residuals $\phi^K=0$ and $\phi^K_\mathrm{GCL}=0 \quad \forall K \in \mathcal{K}^\mathrm{ad}$;
    \item $0.5 < \tau < 1$: all triangles $K \in \mathcal{K}^\mathrm{ad}$ are expanded back from $\mathbf{x}_*$ to the final configuration, and during this deformation  we evaluate the swept areas $\Delta A_i^\mathrm{exp}\quad \forall i \in \mathcal{K}^\mathrm{ad}$ (where the superscript ``exp'' stands for \textit{expansion});
    \item $ \tau = 1$: the final configuration is reached, and we update the areas swept during the fictitious time $\tau=[0,1]$:
    \begin{equation*}
         \Delta A_i \leftarrow \Delta A_i + \Delta A_i^\mathrm{col} + \Delta A_i^\mathrm{exp} \quad \forall i \in \mathcal{K}^\mathrm{ad} \,.
    \end{equation*}
\end{itemize}
This procedure is applied only to the elements involved in the topology modification, that is within the limited region of the grid formed by the union of the elements in $\mathcal{K}^\mathrm{ad}$, as shown in Fig.~\ref{fig:PointAddition},~\ref{fig:PointDeletion} and~\ref{fig:Edgeswapping}.
Indeed, even if the edges external to this region were included, during the collapse and the expansion phases, they would sweep the same area but with opposite sign, so 
$$\Delta A_{ik}^\mathrm{col}+\Delta A_{ik}^\mathrm{exp}=0 \quad \forall e_{ik} \notin \mathcal{K}^\mathrm{ad}\,.$$

Finally, we remark that this procedure is carried out in the same time the mesh adaptation algorithm performs a modification. This can be consider a drawback of our method, but in this way, we can avoid to store the history of the grid modifications and the cumulative property of the swept area makes straightforward the treatment of consecutive topology changes involving the same edges.
Once the mesh adaptation is ended, we compute the interface velocity for all grid nodes through~\eqref{e:defnu_i} considering the whole swept area stored in $\Delta A_i$. Then, we compute the solution on the new grid by solving the non-linear system of equations
\begin{equation}\label{e:finSyst}
    \sum_{K  \mid i \in K} \Phi_i^{\mathrm{HO}, K} = \mathbf{0} \qquad
    \forall i \in \Omega_h^{n+1} \,,
\end{equation}
without any explicit interpolation, but simply considering the new set of DoFs and taking into account the new interface velocities $\nu_i^{n+1}$ while building the space-time residual $\Phi_i^{\mathrm{HO}, K}$.
In this way, the contribution of the topological modifications enters the governing equations as part of the ALE fluxes, as $\nu_i^{n+1}$ are used in~\eqref{FV2RDS_0} to define the  total residuals $\phi^K$.

%% -----------FIGURE 3----------------------
%% -----------END FIGURE 3--------------------

%% -----------FIGURE 4----------------------
%% -----------END FIGURE 4--------------------

%% -----------FIGURE 5----------------------
%% -----------END FIGURE 5--------------------

\subsection{Contribution of modified elements to the discrete governing equations}
To describe how the contributions resulting from the collapse-expansion procedure are included in the governing equations, we take as an example the node insertion, as shown by Fig.~\ref{fig:PointAddition}.
As explained later, the same concepts extend straightforwardly to node deletion and edge swap. 

With reference to Fig.~\ref{fig:PointAddition}, the new node $f$ is inserted by splitting the edge $e_{nj}$, so the two elements sharing this edge form the set $\mathcal{K}^\mathrm{ad}_0=\mathcal{K}^\mathrm{ad}(0\leq \tau< 0.5)$.
Each of these elements is split in two parts, so we have four elements forming the set 
$\mathcal{K}^\mathrm{ad}_1=\mathcal{K}^\mathrm{ad}(0.5 <\tau \leq1)$.
From an implementation point of view, it would be difficult to keep track of the change of connectivity of the elements $K^\prime \in \mathcal{K}^\mathrm{ad}(\tau=0)$, so we introduce the concept of \textit{ghost elements}.
When a triangle undergoes a connectivity change, it is collapsed at $\tau=0.5$, it is marked as ghost element, and it  does not take part in the expansion phase, so it does not belong to the grid $\Omega_h^{n+1}$.
New grid elements $K\in \mathcal{K}^\mathrm{ad}_1$ complying with the new connectivity are generated at $\tau=0.5$, and take part in the expansion phase.
In other words, from the implementation point of view, the triangles forming the set $\mathcal{K}^\mathrm{ad}_0$ are different from the ones forming the set
$\mathcal{K}^\mathrm{ad}_1$, and
$\mathcal{K}^\mathrm{ad}=\mathcal{K}^\mathrm{ad}_0 \cup \mathcal{K}^\mathrm{ad}_1$.
However, since both sets are involved in the collapse-expansion procedure, the elements in both of them generate a contribution to the governing equations:
\begin{itemize}
\item for $K\in \mathcal{K}^\mathrm{ad}_1$, nothing changes from the standard definition of the residual explained in Sec.~\ref{s:ale}, a part that $ \boldsymbol{\eta}_i = 0 \;$ and $|K|=0$ for $t \leq t^{n}$;
\item for $K^\prime \in \mathcal{K}^\mathrm{ad}_0$,
the current integrated nodal normal and the area are null, i.e., 
$\boldsymbol{\eta}_i^{n+1} = 0 \; \forall i \in K^\prime$ and  ${|K^\prime|}^{n+1}=0$, thus the residual of the element $K^\prime$ is (cfr. \eqref{FV2RDS_0} and \eqref{e:timeint})
$$\phi^{K^\prime} = - \sum_{i \in {K^\prime}} \textbf{w}_i  \nu_i \qquad \mathrm{and} \qquad
\Phi_i^{K^\prime} = \frac{{\sum_{q=0}^p} a_q |K^\prime|^{\,n-p} \textbf{w}_i^{n-p}}{\Delta t} + \phi^{K^\prime}_i\,,$$
and, in this case, $\alpha_{\mathrm{LF}} = \max_{i \in K} \vert -\nu_i \vert$. The high order is obtained in the same way as already shown.
Its contribution is purely of ALE nature, thus no stabilization term is required.
\end{itemize}
While assembling the final governing equations for all nodes, in \eqref{e:finSyst}, the summation is performed by considering the nodal residuals $\Phi_i^K$ of all elements that are \textit{or were} connected to $i$.

Care should be taken during the time integration, as the number of DoFs changes between steps: the system of governing equations~\eqref{e:finSyst} is augmented when a new node is inserted. For the new equation, the time integration scheme in~\eqref{e:timeint}, as well as~\eqref{e:dgcl}, does not include the contributions of the previous steps because they are multiplied by the area $|K|$ of elements that were not present in the grid for $t<t^{n+1}$, so their area is null.
On the other hand, the ghost elements generated during the collapse phase must be kept in the system of governing equations for the next $p$ steps.
To clarify, consider to compute the solution at $t^{n+1}$ with the BDF scheme in~\eqref{e:timeint} with two steps ($p=1$).
Ghost elements $K^\prime$ removed during the previous time step (from $t^{n-1}$ to $t^n$) have a non-null area ${|K^\prime|}^{n-1}>0$ and null areas  $|K^\prime|^{n}={|K^\prime|}^{n+1}=0$, so their contribution to the residuals $\Phi_i^{K^\prime}$ is not null because of the interface velocities resulting from the area swept during the collapse phase ($\Delta A^n\neq 0$, but $\Delta A^{n+1}=0$). From the next step, i.e., while computing the solution at $t^{n+2}$, these ghost elements will not contribute anymore to the governing equations, so they can be completely removed.

\paragraph{Extension to node deletion}
The same considerations about the treatment of the connectivity change and the residual contribution for the case of node insertion can be done also for node deletion.
With reference to Fig.~\ref{fig:PointDeletion}, the set $\mathcal{K}^\mathrm{ad}_0$ contains all triangles sharing the node $f$, the one to be deleted.
At $\tau=0.5$, two triangles and the node $f$ are deleted. From the implementation point of view, all triangles in $\mathcal{K}^\mathrm{ad}_0$ are marked as ghost elements and are not expanded, and new triangles conforming with the new topology are created and form the set $\mathcal{K}^\mathrm{ad}_1$.
The total residual of the elements of both $\mathcal{K}^\mathrm{ad}_0$ and $\mathcal{K}^\mathrm{ad}_1$ are computed taking into consideration the observations done before.
In addition, the solution over the deleted node should be computed because of conservativity constraints, since it may contribute to some ghost elements that contribute to a node $i$ still present in the grid. The removed nodes are therefore considered as \textit{ghost nodes} for the next $p$ steps after their deletion. Afterwards, they are actually deleted from the solution procedure.

\paragraph{Extension to edge swap}
The edge swap is described in a very similar way:  as shown in Fig.~\ref{fig:Edgeswapping},
the set $\mathcal{K}^\mathrm{ad}_0$ contains the two triangles sharing the edge $e_{nj}$, at $\tau=0.5$ the edge is deleted and a new edge $e_{im}$ is generated to connect the other two nodes in the quadrilateral formed by the union of the two triangles.
From the implementation point of view, the triangles initially sharing the edge $e_{nj}$ become ghost elements, while two new ones are generated to comply with the new connectivity, and form the set $\mathcal{K}^\mathrm{ad}_1$.
This case is somehow simpler than the previous ones as no nodes are added or deleted, while the considerations about ghost and new elements done for the node insertion apply here, too.
Then, their contributions to the nodal residual are computed in the same way.

\section{Summary of the computational procedure}\label{ss:compProc}
Before presenting the numerical results that show the validity of the proposed approach, we give some information about the computational workflow and the used software. 

As already anticipated, in addition to time-dependent simulations,
we present also two steady problems.
In adaptive steady simulations, $\mathcal{M}$ cycles \textit{mesh adaptation--solution computation} are performed to iteratively optimize the mesh according to the steady solution which is re-computed at every step over the new grid starting from the old one. To map the solution between different grids without interpolation, we use at each step the ALE RD scheme described in the previous section, considering an artificial time step $\Delta t =1$ and an implicit Euler scheme as explained in Sec.~\ref{ss:timesteady}.

In unsteady simulations, the time has the standard physical meaning and we reach the final solution $\mathit{w}(\mathbf{x},t_\mathrm{F})$ starting from the initial solution  $\mathit{w}(\mathbf{x},t_\mathrm{0})$ according to a time-marching process, splitting the time interval in $N_t$ steps, so that $\Delta t = (t_\mathrm{F} - t_\mathrm{0}) /N_t$.
In addition, at each time step, to solve system \eqref{e:finalsystemRD}, we use an implicit pseudo-time stepping technique, using the same BDF scheme adopted for the physical time integration, but with a variable integration step $\Delta \tau$, computed according to an adaptive CFL law~\cite{Re2016}.
In unsteady simulations, we can choose to perform an adaptation step at each physical time step or every a certain number.

Each adaptation step can be summarized as follows. 
\begin{enumerate}
    \item Starting from an initial mesh $\Omega_h^n$ and the associated solution $\mathit{w}^n=\mathit{w}(\mathbf{x},\Omega_h^n)$, the target grid spacing $\chi_h^{n+1}$ is computed according to the derivatives of the flow variables.
    \item The re-mesher library is called, with the mesh  $\Omega_h^n$ and the target grid spacing $\chi_h^{n+1}$ as inputs. Here, an inner loop is performed:
    \begin{itemize}
        \item the re-mesher cycles over the grid entities and checks if the connectivity should be modified to comply with the target spacing,
        \item when the re-mesher inserts or removes a node or swaps an edge, it communicates to the flow solver the information about the topology modification and stays idle,
        \item the flow solver computes the swept areas executing the collapse-expansion procedure,
        \item the re-mesher resumes from idle and checks next entity.
    \end{itemize}
    \item The interface velocities $\nu_i^{n+1}$ over the adapted grid $\Omega_h^{n+1}$ are computed.
    \item The solution $\mathit{w}^{n+1}$ of the Euler equations is evolved over the new grid $\Omega_h^{n+1}$  using the BDF ALE RD scheme described in Sec.~\ref{s:ale}.
    \item The stored results are updated (remind that, using the second-order BDF scheme, also $\mathit{w}^{n-1}$ and $\Omega_h^{n-1}$ are needed to compute $\mathit{w}^{n+1}$).
\end{enumerate}
This procedure is implemented in the in-house FORTRAN solver,  \textsf{Flowmesh}~\cite{Re2017aiaa,Re2016,IsolaPhd}, 
which is linked to the re-mesher library \textsf{mmg}~\cite{Dobrzynski2008,Dapogny2014,mmgWeb}.

\section{Benchmark tests considering analytical solutions}\label{s:validation}
In this section, we present the numerical results of three tests (two unsteady and one steady), for which the analytical solution is known.
The initial data and the quantitative results presented in this section as dimensionless are scaled with respect to the mean sea level condition of ICAO standard atmosphere.

\subsection{Numerical validation of Geometric Conservation Law}\label{ss:valid_gcl}
\begin{figure*}
\centering
\includegraphics[width=\textwidth]{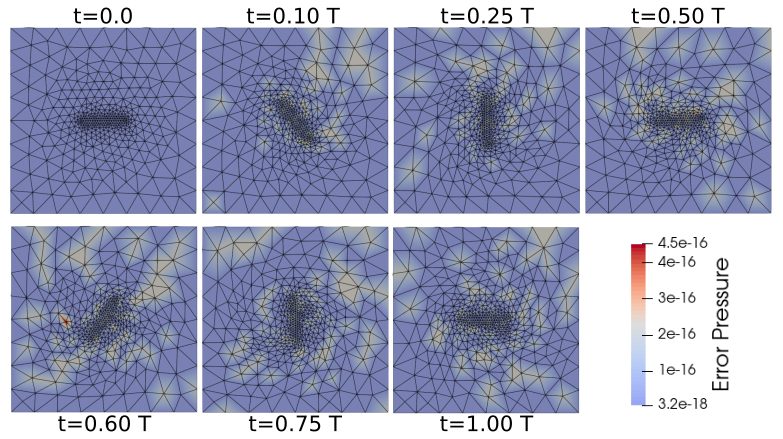}
\caption{Numerical validation for GCL in the adaptive simulation: 
computational grid at different times. The domain is colored according to the nodal error on the pressure, $Err_P(\mathbf{x}_j, t^n) = \vert P_j(t^n) - P_0\vert$.}
\label{f:gcl_grid}
\end{figure*}

The fulfilment of the GCL can be checked numerically by computing the evolution of
a uniform field over a dynamic grid.
For this purpose, we consider the square computational grid $\Omega_h^0= [-1 \,, 1]^2$ split in two regions (see Fig.~\ref{f:gcl_grid}): a box corresponding to the central rectangular area with corners $[-0.25 \,, -0.025]$ and $[0.25 \,, 0.025]$, where the grid elements have a characteristic edge length $h_0^\mathrm{box}=0.025$; and the surrounding part, where the element size progressively increases from $h_0^\mathrm{box}$ to $h_0^\mathrm{max}=0.25$, imposed on the grid boundary $\partial\Omega_h^0$.
A rigid pitching motion is imposed to the box region:
$\alpha(t) = A \sin(\omega t)$, with $A=90^\circ$ and $\omega = 2 \pi / 100$, while the grid boundary $\partial\Omega_h^0$ is kept fixed. The grid elements in the surrounding part are deformed and/or adapted to account for the box motion, preserving the initial grid spacing.
We initialize the domain with uniform density $\rho_0=1.0$, velocity $\mathbf{u}_0 = [1.0 \,, 1.0 ]$ and pressure $P_0=1.0$, and compute the evolution of the solution for a pitching time period. The final time is reached in 120 steps, corresponding to a CFL number of $11.8$ with respect to $h_0^\mathrm{box}$.

\begin{figure}
\centering
\includegraphics[width=0.5\textwidth]{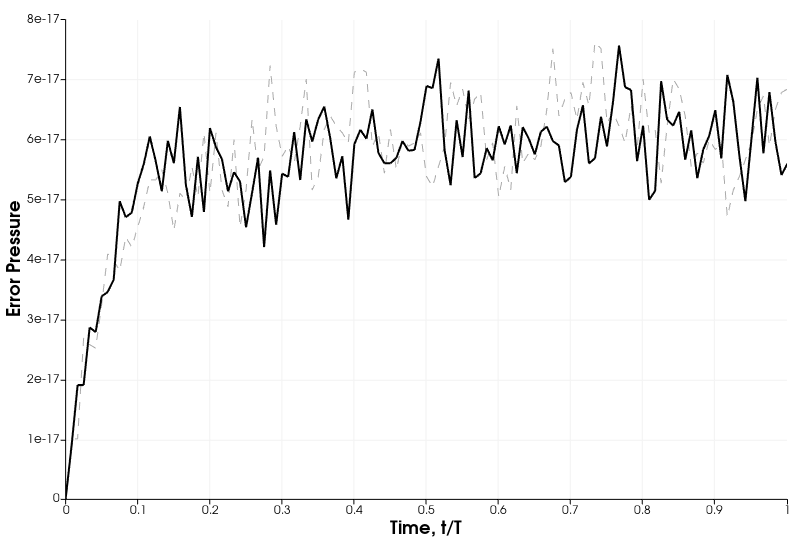}
\caption{Numerical validation for GCL: error evolution during the simulation.
The plot shows the error on the pressure, averaged over all grid points, in the adaptive simulation (thick, solid line) and in a simulation where only mesh deformation is used (grey, dashed line).
The maximum values of the nodal error measured in the simulations are both around $4.4 \,\cdot 10^{-16}$.}
\label{f:gcl_error}
\end{figure}

Two simulations are performed: a first-one where the grid connectivity is kept fixed, but the grid is deformed according to the elastic analogy to cope with the rigid movement of the elements in the box region; a second one where the grid is adapted every 10 time steps aiming for an improved grid quality preserving the actual grid spacing.
The error of the final solution with respect to the initial uniform flow field, shown in Fig.~\ref{f:gcl_error}, is of the order of $10^{-16}$, close to machine precision.

%% -----------FIGURE 6----------------------
%% -----------END FIGURE 6------------------

%% -----------FIGURE 7----------------------
%% -----------END FIGURE 7------------------

\subsection{Isentropic vortex}
%% -----------FIGURE 8----------------------
\begin{figure}
\centering
\includegraphics[width=0.5\textwidth]{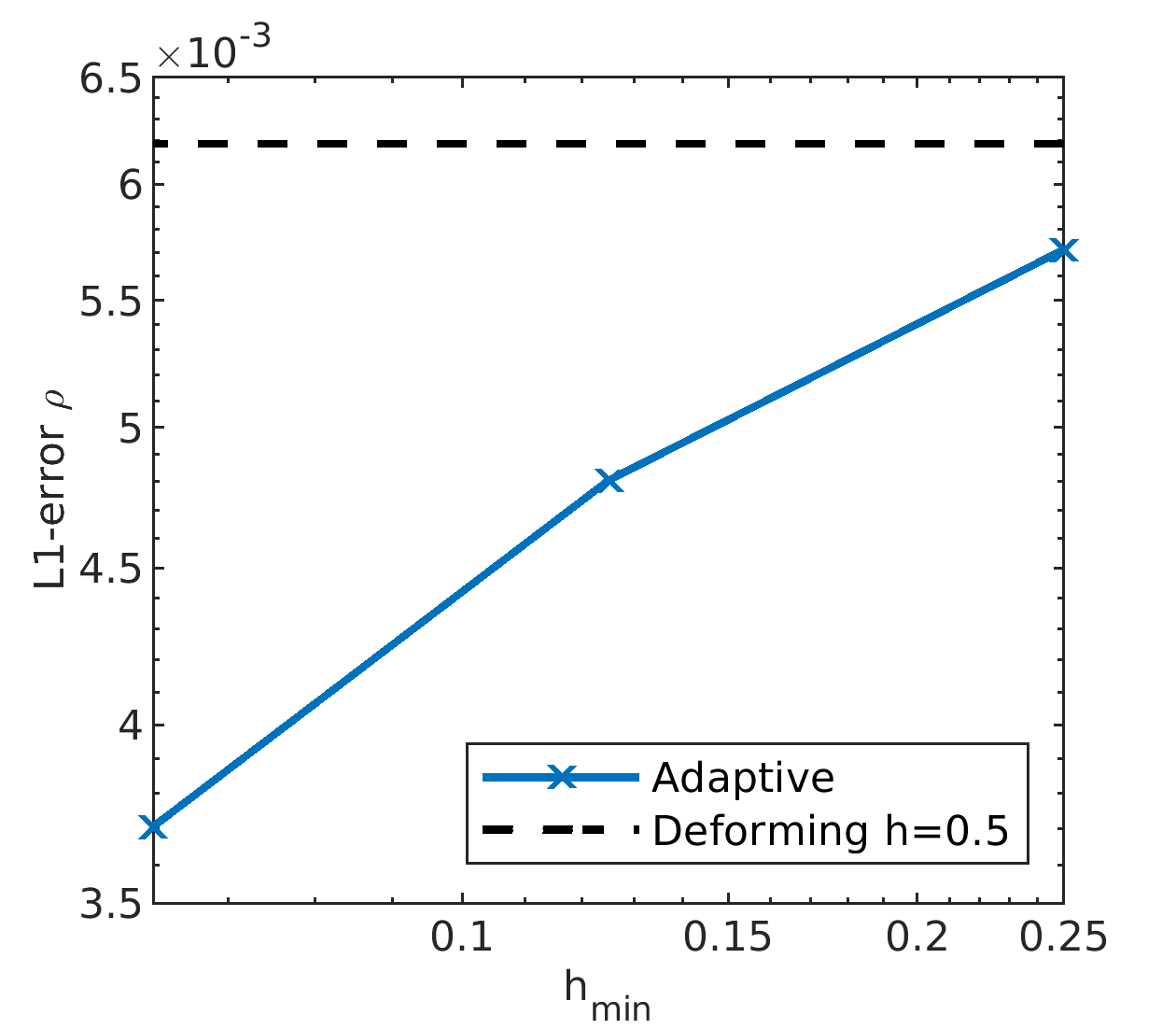}
\caption{Isentropic vortex test: $L^1$ norm of the error on the density at the final time ($t^f=10)$ for different minimal edge sizes.
The starting grid has a uniform grid spacing $h_0 = 0.5$. The adaptive simulations are carried out imposing minimal edge sizes of $h=[0.25, 0.125, 0.0625]$, keeping constant the CFL number with respect to this size.
For a comparison, the dashed line shows the error obtained over the initial grid deforming at fixed connectivity and preserving the initial grid spacing.}
\label{f:vortex_err}
\end{figure}
%% -----------END FIGURE 8------------------

The isentropic vortex problem is a benchmark test for which the analytical solution is known and regular~\cite{Shu1998}.
The initial solution is defined as
\[ \rho_0=\rho_\infty + \delta \rho, \quad
\mathbf{u}_0 = \mathbf{u}_\infty + \delta \mathbf{u},\quad
P_0 =P_\infty + \delta P,
\]
where the undisturbed field is characterized by density $\rho_\infty = 1.0$, velocity $\mathbf{u}_\infty = [1.0 \,, 1.0 ]$, pressure $P_\infty=1.0$, and temperature $\theta_\infty=1.0$.
The vortex is centered at $[x_\mathrm{v}, y_\mathrm{v}]$ and the perturbations are defined as
\[ 
\delta \theta = -\frac{(\gamma-1) \lambda^2}{8 \gamma \pi^2}\exp{\left(1-r^2\right)}\,,
\qquad
\delta\mathbf{u} =\begin{bmatrix}
\delta u\\ \delta v
\end{bmatrix} = \frac{\lambda}{2 \pi} \exp{\left(\frac{1-r^2}{2}\right)} \begin{bmatrix}
-(y-y_\mathrm{v})\\x-x_\mathrm{v}
\end{bmatrix} \,,\]
where $r= \sqrt{(x-x_\mathrm{v})^2+(y-y_\mathrm{v})^2}$, and $\lambda=5.0$ is the vortex strength.
The perturbation on the density and on the pressure are defined as
$ \delta \rho = \left( 1 + \delta \theta \right)^{\frac{1}{\gamma -1}} -1$, 
and $ \delta P = \left( 1 + \delta \theta \right)^{\frac{\gamma}{\gamma -1}} -1$.

To run this test, we use a computational grid $[-10,10]^2$ and the initial position of the vortex is $x_\mathrm{v}^0=-5$, $y_\mathrm{v}^0=-5$, and we compute the evolution until the final time $T^\mathrm{f}=10$. The analytical solution consists in the advection of the initial vortex until its center is $x_\mathrm{v}^\mathrm{f}=5$, $y_\mathrm{v}^\mathrm{f}=5$.
At the final step, we compute the $L^1$ norm of the error on the density, as
$$ \mathrm{Err}_\rho = 
\frac{\int_{\Omega^h} \vert \rho^h(\mathbf{x}) - \rho^\mathrm{an}(\mathbf{x}) \vert \, \mathrm{d}\mathbf{x}}
{ \int_{\Omega^h} \rho^\mathrm{an}(\mathbf{x})\, \mathrm{d}\mathbf{x} } \,, $$
where $\rho^\mathrm{an}$ is the analytical solution.

To show the effect of the adaptation over the error, we started from a uniform grid, with $h_0=0.5$, and we run three adaptive simulations, imposing a different minimal edge size. The imposed minimal edge size is reached in the vortex region, while far from it the initial grid spacing is preserved. Figure~\ref{f:vortex_err} shows how the error decreases with the minimal edge size. A proper convergence analysis is not performed, as the final grids do not present a uniform grid spacing, and they are not refined in a hierarchical way.

\subsection{Oblique shock reflection}\label{ss:schok}

%% -----------FIGURE 6----------------------
\begin{figure*}
\centering
\includegraphics[width=0.95\textwidth]{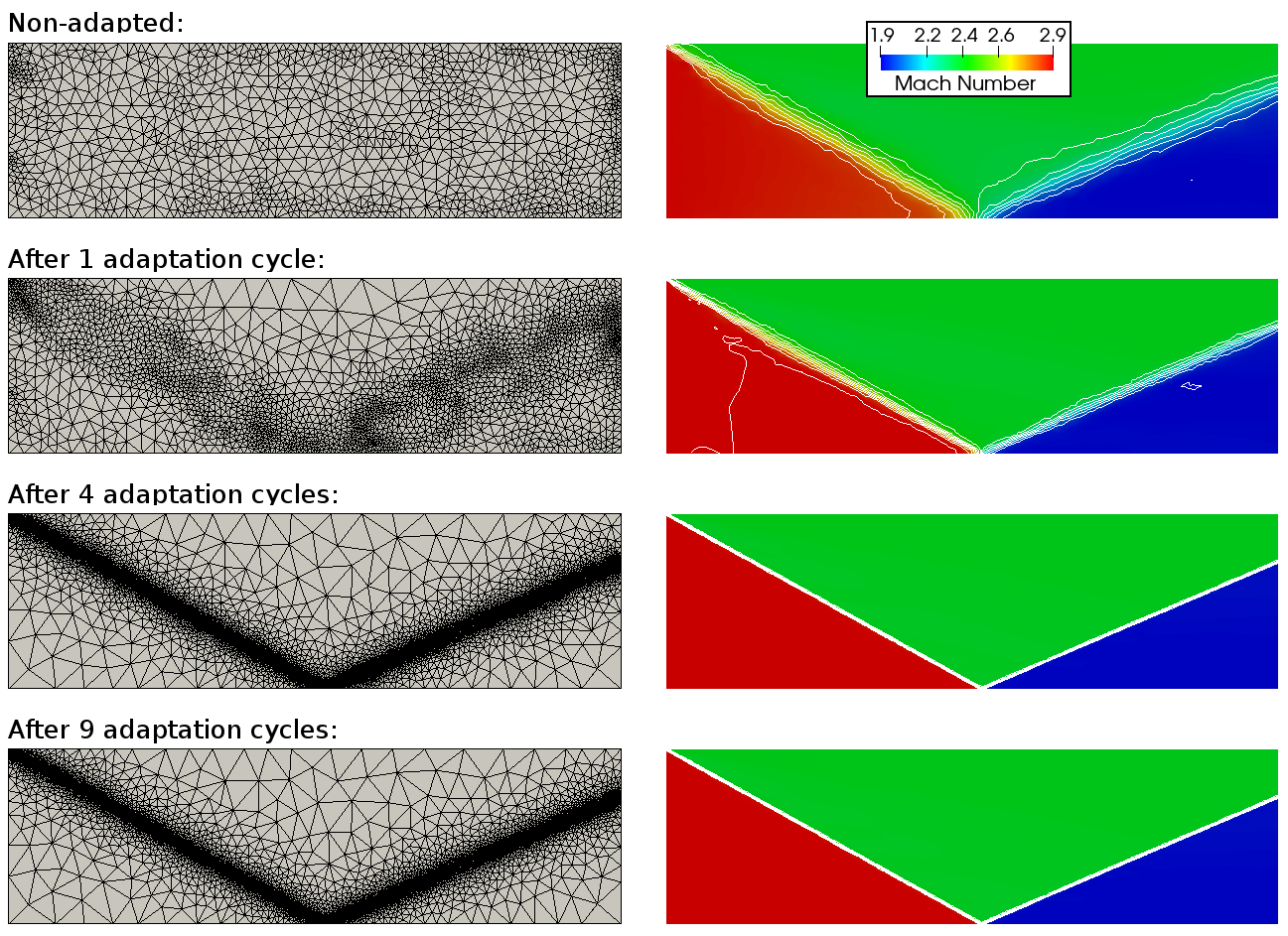}
\caption{Grid and Mach number contour plots of the simulation of the oblique shock reflection. 
The first row shows the initial, non-adapted grid, composed of 615 nodes, and the associated solution.
The next rows show the grid and the solution after 1, 4, and 9 adaptation cycles.
The number of nodes are, respectively, 1230, 8239 and 7213.}
\label{f:shock1}
\end{figure*}
%% -----------END FIGURE 6------------------

%% -----------FIGURE 7----------------------
\begin{figure*}
\centering
\includegraphics[width=\textwidth]{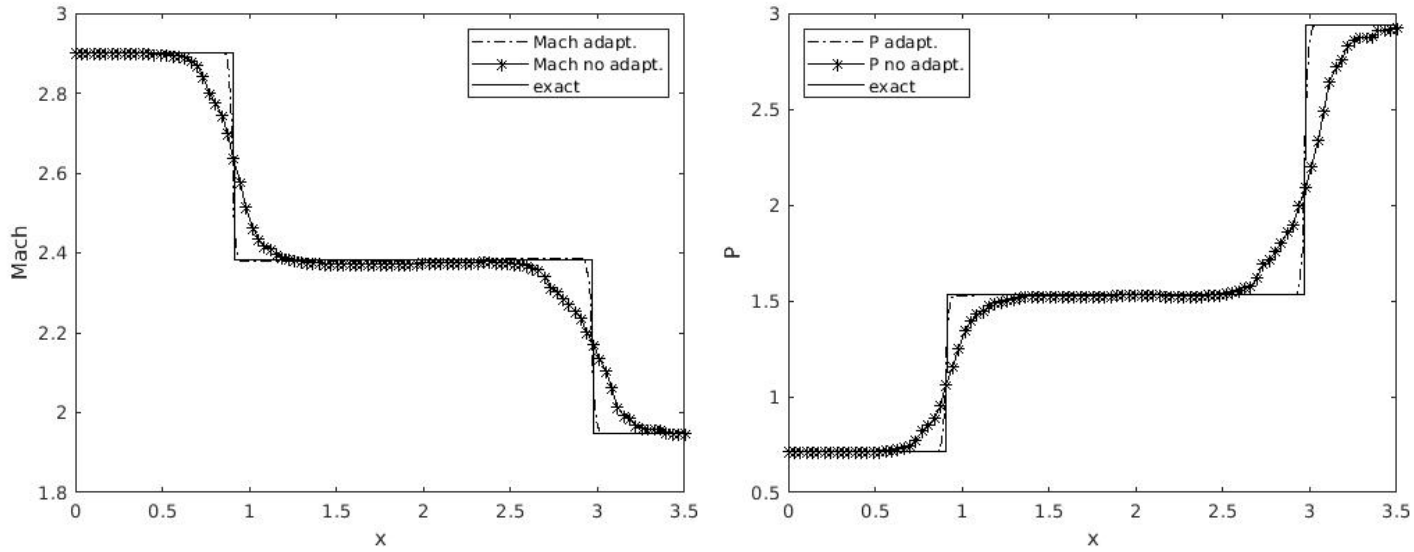} \\
(a) Mach number \hspace{5cm} (b) Pressure
\caption{Results of the simulation of the oblique shock reflection, extracted along the horizontal line at $y=0.5$.
The Mach number and the pressure along $x$ for the solution computed on the initial grid and on the grid after 9 adaptation cycles are compared to the analytical one, computed according to the Rankine-Hugoniot relations.}
\label{f:shock2}
\end{figure*}
%% -----------END FIGURE 7--------------------
%% -----------FIGURE 8----------------------
\begin{figure*}
    \centering
    \includegraphics[width=\textwidth]{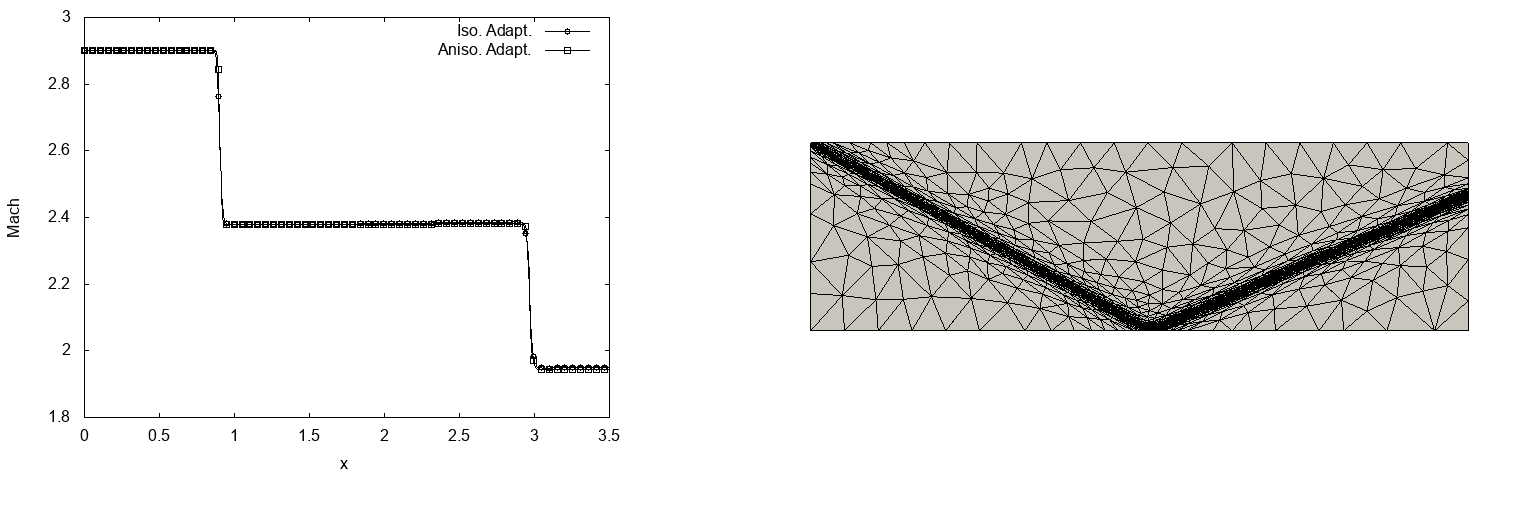} \\
    \hspace{-1cm} (a) Mach number \hspace{6.5cm} (b) Grid
    \caption{Results of the simulation of the oblique shock reflection using anisotropic mesh adaptation. (a) Mach number profile extracted along the horizontal line at $y = 0.5$, compared to the one obtained with isotropic  adaptation. (b) Grid resulting after 9 cycles of anisotropic adaptation. The number of nodes is 2528.}
    \label{f:shock3}
\end{figure*}
%% -----------END FIGURE 8--------------------

The last benchmark test concerns the so-called regular shock reflection, which involves two shock waves: an incoming shock wave incident a planar solid wall and an oblique reflected wave, which intersect at a point over the plane surface.
In this steady test, mesh adaptation is important to resolve accurately the shock, without knowing a-priori its position. Indeed, starting from a uniform initial solution, the shock position and resolution change significantly, especially during the first adaptation steps.

The domain is rectangular, the flow enters from the left, the wall is at the bottom and the right side is a free outlet.
The initial conditions, uniform, are  $\rho_0=1.6997$, $\mathbf{u}_0 = [2.619 \,, -0.5069 ]$ (corresponding to a Mach $M_0=2.38$), and $P_0=1.529$.
The same pressure $P_0$ is impose on the top, while on the left boundary we impose $\rho_\mathrm{IN}=1.0$, $P_\mathrm{IN}=0.714$, and a horizontal flow at Mach number $M_\mathrm{IN}=2.9$.

We execute this simulation performing $\mathcal{M}=9$ cycles of adaptation, beyond which the solution does not change significantly. 
The results of the simulation at different steps are shown in Fig.~\ref{f:shock1}, where we can appreciate the positive effects of mesh adaptation on the shock resolution. The adaptation criterion is based on the Hessian matrix of the pressure.

For this simple fluid problem, under the hypotheses valid for the Euler equations, we can compute the analytical solution using the Rankine-Hugoniot conditions on shock front (e.g., see~\cite{Thompson1972}).
Hence, it is possible to compare the result  along $x$ for a fixed $y$ coordinate. Figure~\ref{f:shock2} shows that the solution on the adapted grid reaches an excellent agreement with the analytical one.

We run also an anisotropic adaptation, using as error indicator the Hessian of the pressure field. The starting grid, the simulation parameters and the number of adaptation cycles are exactly the same as the ones described above. The results of this simulation are compared to the results of the isotropic adaptive simulation in Fig.~\ref{f:shock3}. The agreement with the analytical solution is satisfactory for both cases, but what changes is the number of nodes after the adaptation: after 9 cycles, the number of nodes reached with isotropic adaptation is 7213, whereas in the anisotropic case it is equal to 2528, almost one third the former one. This test shows that the proposed method can be used straightforwardly also with anisotropic adaptation, which, in general, provides a saving of computational time, as it can reach the same level of accuracy of isotropic adaptation using less nodes. 

\section{Numerical results}\label{s:results}
This section illustrates the results of three adaptive simulations, to support the validity of the proposed methods and its applicability to unsteady compressible flow fields.
As in Sec.~\ref{s:validation}, the dimensionless data are scaled with respect to the mean sea level condition of ICAO standard atmosphere.

\subsection{Steady flow over the RAE2822 airfoil}\label{ss:rae}
\begin{figure*}
\centering
\includegraphics[width=0.45\textwidth]{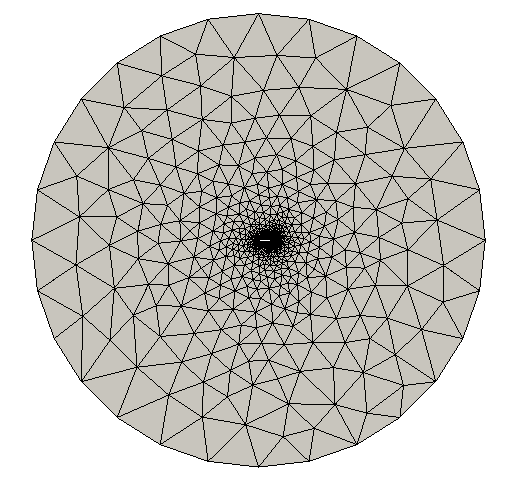} \hfill
\includegraphics[width=0.5\textwidth]{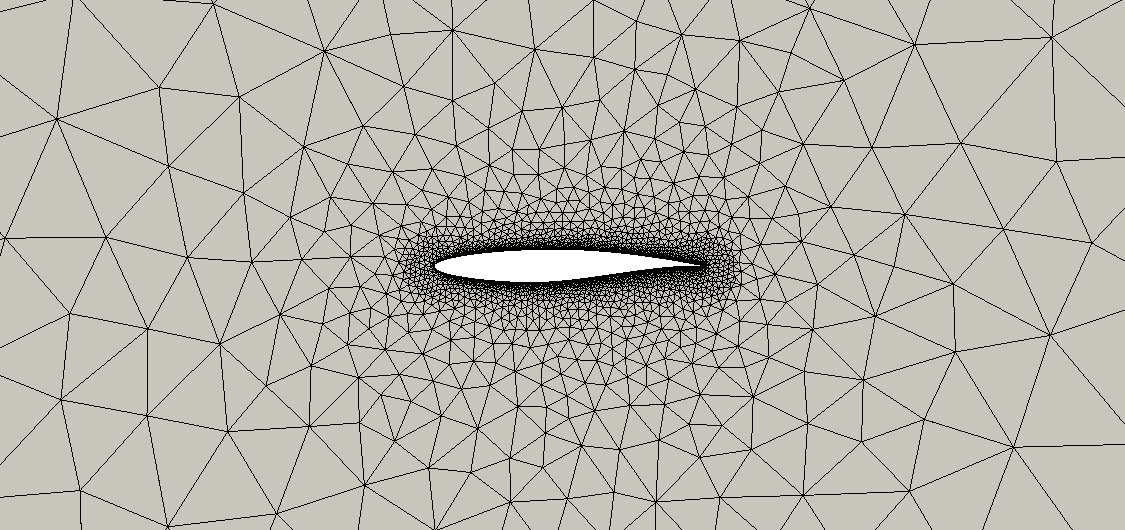}
\caption{Initial mesh for the RAE2822 test case: on the left there is the whole grid composed of 2713 nodes, on the right a detail of the region close to the airfoil.}
\label{f:raemesh0}
\end{figure*}

This test involves the transonic flow over the RAE-2822 airfoil.
The computational domain is an O-grid with a radius of $12$ times the chord $c$ of the airfoil.
The initial mesh, shown in Fig.~\ref{f:raemesh0}, is composed of 2713 nodes.
The free-stream conditions, namely density $\rho_\infty=1.4857$, Mach number $M_\infty=0.729$, and angle of attack $\alpha=2.31^\circ$, are imposed on the external boundary, while the airfoil is modeled as an inviscid wall.

\begin{figure*}
\centering
\includegraphics[width=0.7\textwidth]{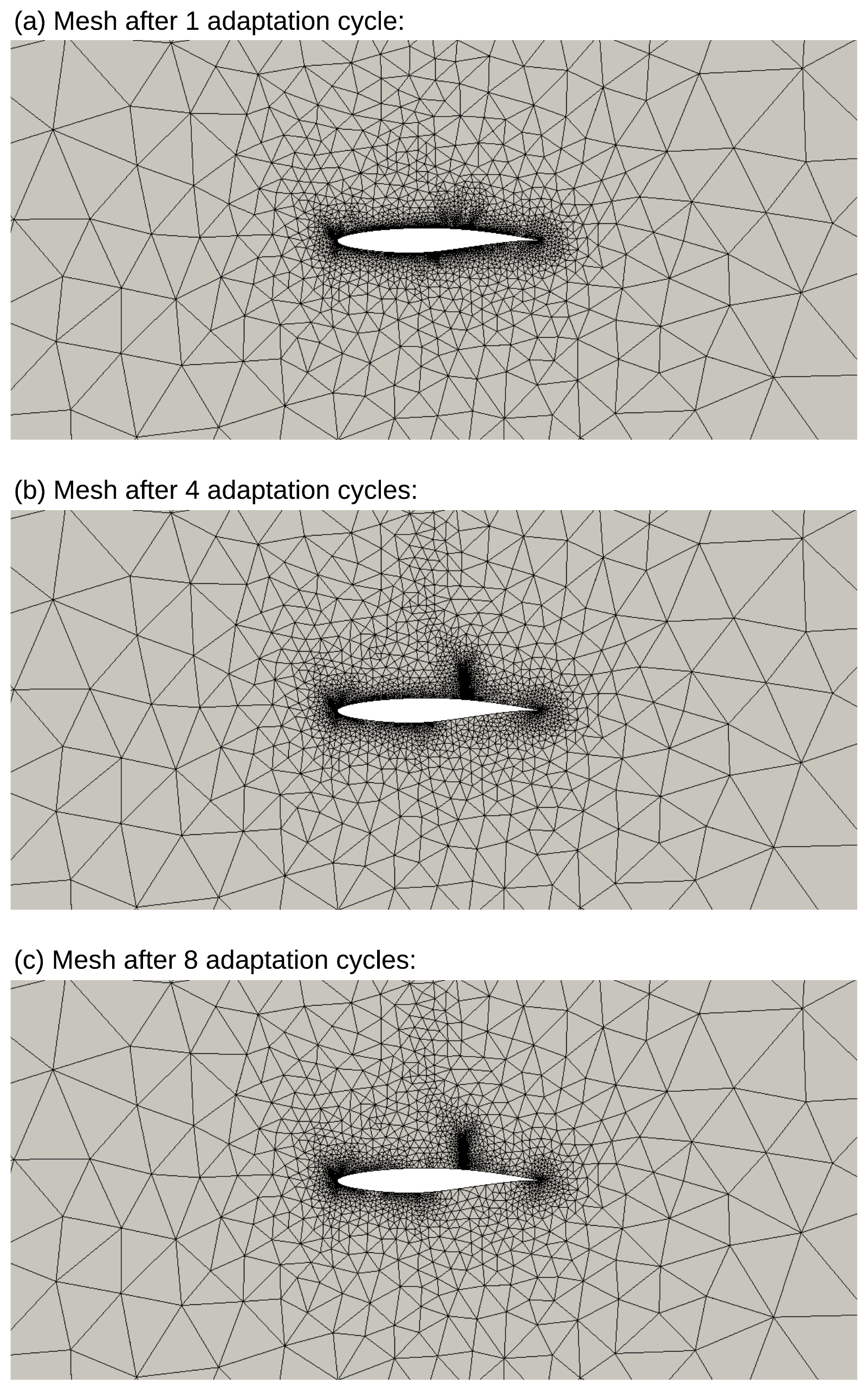}
\caption{Adapted meshes during the RAE2822 test case: detail of the region close to the airfoil.
From top to bottom, the grid are composed of 3402, 3841 and 3359 nodes.}
\label{f:raemesh}
\end{figure*}

\begin{figure*}
\centering
\includegraphics[width=0.7\textwidth, height=7cm]{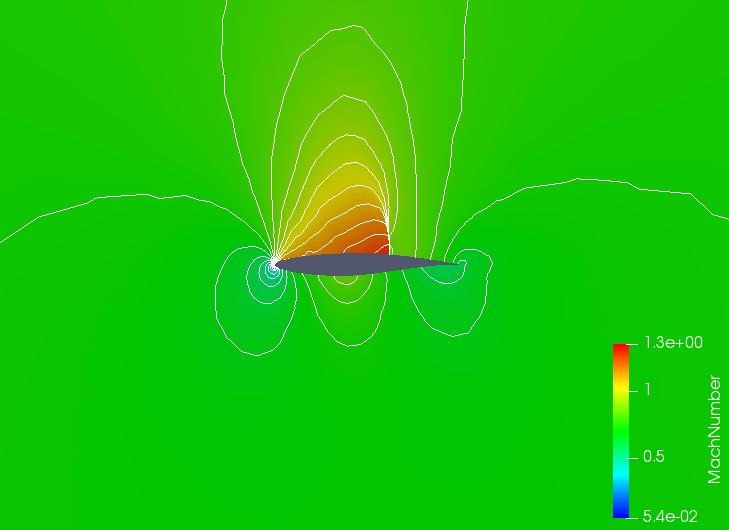} 
\caption{Solution of the simulation of the RAE2822 test case: contour plot of the Mach number along with the Mach-isoline in the domain region close to the airfoil after $\mathcal{M}=8$ cycles of \textit{mesh adaptation--solution computation}.
The related mesh is the one shown in the last row of Fig.~\ref{f:raemesh}: we can observe the mesh refinement in correspondence of the shock wave.}
\label{f:raesol}
\end{figure*}

We compute the solution performing $\mathcal{M}=8$ adaptation cycles, using as adaptation criterion the Hessian matrix of the pressure.
The meshes after different adaptation steps are shown in Fig.~\ref{f:raemesh}, while the flow field computed on the final mesh is shown in Fig.~\ref{f:raesol}. Due to the transonic condition, a shock wave forms on the upper side of the airfoil. The mesh is refined in accordance with this behavior, preserving also a small grid spacing close to the leading and trailing edge, where the curvature of the body is maximum.
From Fig.~\ref{f:raemesh}, we can appreciate also the effect of mesh coarsening, especially in the lower back part of airfoil, where the flow field does not  exhibit significant gradients (see Fig.~\ref{f:raesol}). Indeed, the final mesh contains almost the same number of nodes as the one after the first adaptation step, but they are distributed in a more effective way, according to the flow features.

\begin{figure*}
    \centering
    \includegraphics[width=0.7\textwidth,height=7cm]{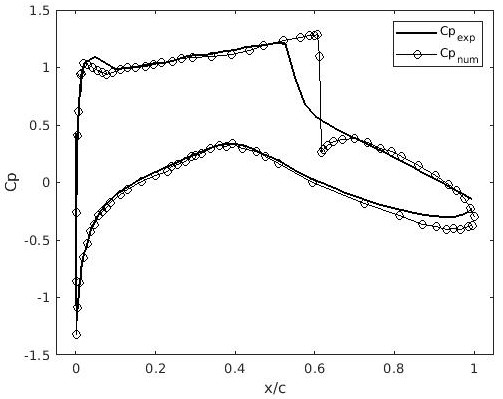}
    \caption{Pressure coefficient over the RAE2822 airfoil: the numerical value computed in the inviscid adaptive simulation is compared to the experimental data given in~\cite{Cook1979}.}
    \label{f:raecp}
\end{figure*}

Finally, we compare the numerical solution with the experimental data provided in~\cite{Cook1979}.
In particular, we compare the value of the pressure coefficient $C_P$ over the airfoil, computed as:
$$ C_P(\mathbf{x}_i) = \frac{P_i - P_\infty}{0.5 \rho_\infty \Vert \mathbf{u}_\infty \Vert^2} , 
\qquad  \forall \mathbf{x}_i \in \partial \Omega_h^\mathrm{b} $$
where $\Omega_h^\mathrm{b} $ is the numerical boundary representing the airfoil.
This comparison is shown in Fig.~\ref{f:raecp}. There is a good agreement, taking into consideration the inviscid hypotheses behind the Euler equations.

%% -----------FIGURE 9----------------------
%% -----------END FIGURE 9------------------

%% -----------FIGURE 10----------------------
%% -----------END FIGURE 10------------------

%% -----------FIGURE 11---------------------
%% -----------END FIGURE 11-----------------

%% -----------FIGURE 12---------------------
%% -----------END FIGURE 12-----------------

\subsection{Mach 3 forward facing step}\label{ss:step}
The channel forward facing step~\cite{Woodward1984} test case is used to asses the capability of the proposed ALE RD scheme to capture complex structures arising from the interaction between shock waves. 
The Mach number for this test is set equal to 3, corresponding to the inlet flow imposed on the left boundary. The upper and lower boundaries are modelled as inviscid walls.

\begin{figure*}
    \centering
    \includegraphics[scale=0.7]{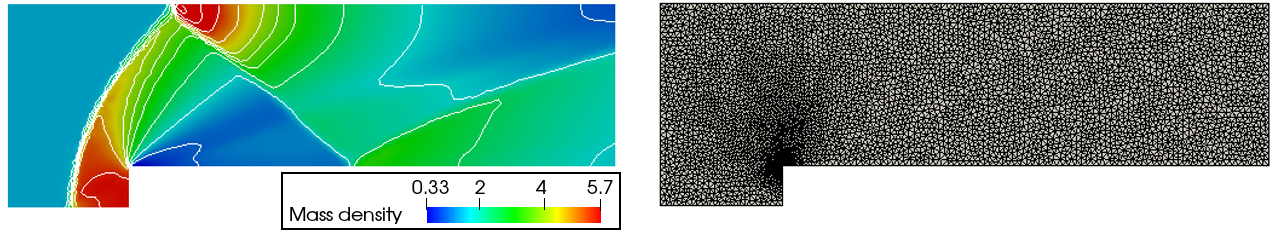}
    \caption{Density contours for a simulation performed without mesh adaptation, at time step $t= 2.31$, over the fixed mesh shown on the right (made of 4668 nodes), which is used also as initial mesh for the adaptive simulation.}
    \label{fig:StepNoAdapt}
\end{figure*}

We investigate the transient evolution of the flow, solving the unsteady Euler equations up to a final time of $t_F=3.84$. At each time step, an implicit pseudo-time stepping technique with adaptive CFL is used to solve the non-linear systems of equations.
The starting CFL number used for the implicit solver is set equal to 0.75 and the adimensional time step is set as $\Delta t = 0.03$. 
Figure~\ref{fig:StepNoAdapt} displays the initial grid, made of 4668 nodes, and the results of a non-adaptive simulation.

\begin{figure*}
\centering
\includegraphics[scale=0.7]{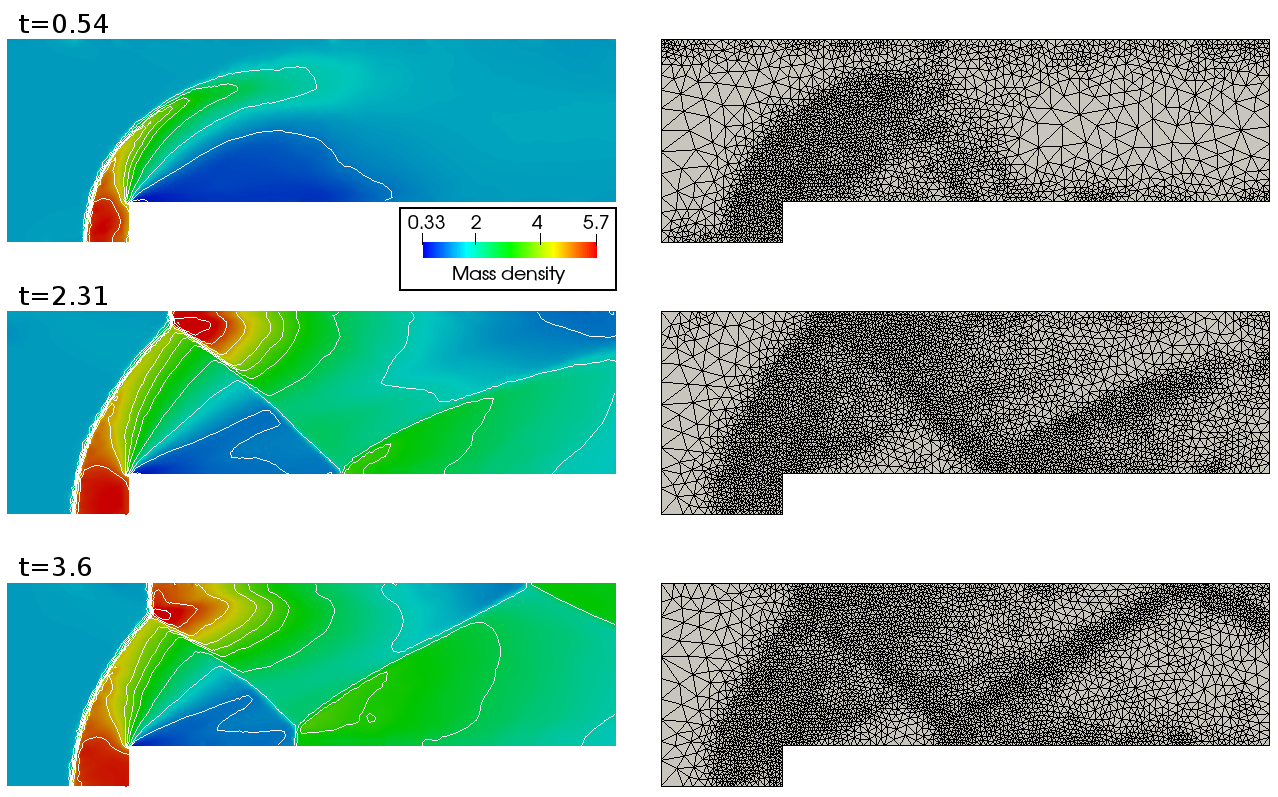} 
\caption{Adaptive simulation of the forward facing step at different time steps: from top to bottom, $t = 0.54$, $t = 2.31$, and $t = 3.6$. On the left side, density contour plots; on the right side, the corresponding grids, made of, respectively, 2843, 4338 and 4385 nodes.}
\label{f:AllSteps}
\end{figure*}

The results of the adaptive simulation at different time steps are shown in
Figure~\ref{f:AllSteps}.
Mesh is adapted at each time iteration, leading to a refinement close to the regions where shocks occur and close to the sharp corner. 
The fluid structures are well captured and we can also notice that on the left hand side of the channel the mesh is not refined due to the absence of strong gradients.
Moreover, comparing the density fields in the last row of Fig.~\ref{f:AllSteps} and in Fig.~\ref{fig:StepNoAdapt}, the adaptive simulation is able to accurately reproduce the flow features using a lower number of grid points.

%% -----------FIGURE 14---------------------
%% -----------END FIGURE 14-----------------
%% -----------FIGURE 13---------------------
%% -----------END FIGURE 13-----------------

\subsection{Unsteady flow around the NACA0012 pitching airfoil}\label{ss:pitch}
A test involving boundary displacement is now considered, with the simulation of the flow around the oscillating NACA0012 airfoil, at a free-stream Mach ${M} = 0.755$. 
The amplitude of the pitching motion is described by the equation $\alpha = \alpha_{\infty} + \alpha_0 \sin{(\omega t)}$, with $\alpha_{\infty} = 0.016^{\circ}$ and $\alpha_0=2.51^\circ$. The final simulation time is $t_F = \frac{4 \pi}{\omega} = 87.1778$ and the number of time steps is $200$. The initial grid used for this simulation has the same characteristics of the one reported for the test case \ref{ss:rae}. The Hessian of the Mach number is the parameter used as mesh adaptation criterion.

\begin{figure*}
\centering
\includegraphics[width=0.45\textwidth]{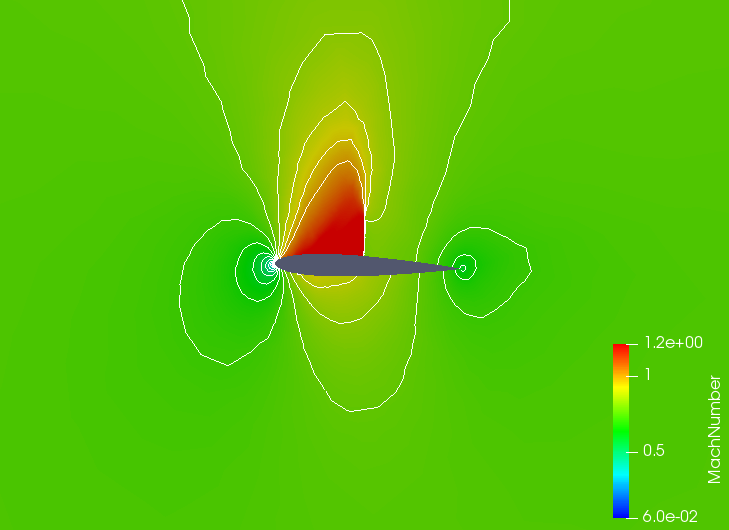} \hspace{0.5cm}
\includegraphics[width=0.45\textwidth]{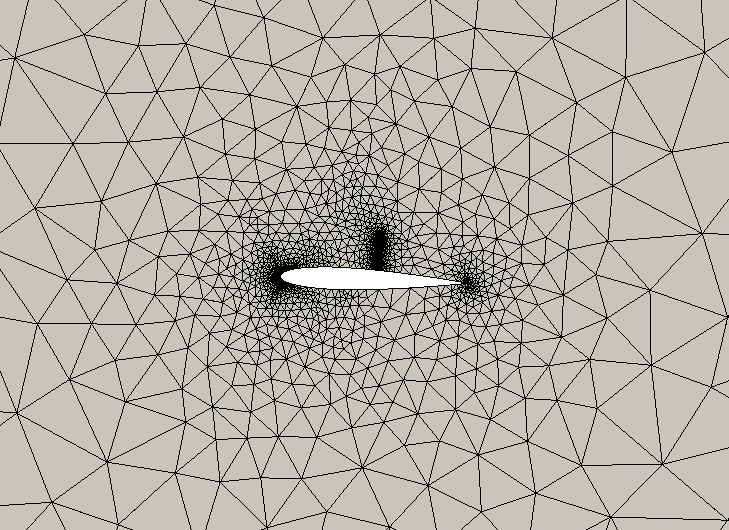}\\[0.5cm]
\includegraphics[width=0.45\textwidth]{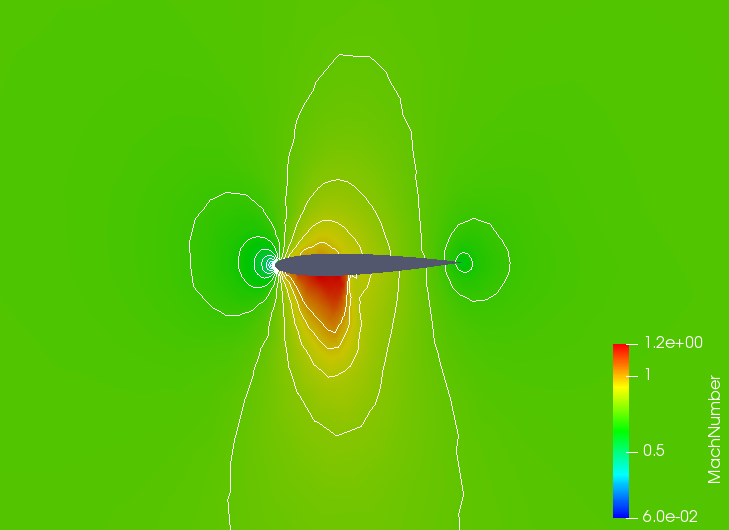} \hspace{0.5cm}
\includegraphics[width=0.45\textwidth]{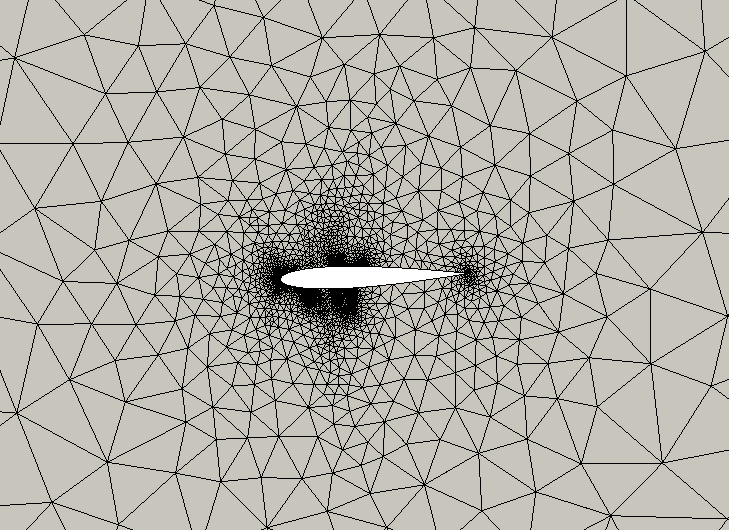}
\caption{Results of the NACA0012 pitching airfoil: Mach contour plots and grids at the maximum pitching angle (time $t=\frac{5}{2}\frac{\pi}{\omega}$) in the first row, and at the minimum pitching angle (time $t=\frac{7}{2}\frac{\pi}{\omega}$) in the second row. The grids are composed respectively by 2508 and 3174 nodes respectively. The starting grid has 1537 nodes.}
\label{f:nacaMINMAX}
\end{figure*}

Figures~\ref{f:nacaMINMAX} shows the results at the maximum and minimum pitching angle, along with the corresponding grid. The proposed ALE RD scheme is able to capture the shocks arising on the top and the bottom of the airfoil. As expected, the upper shock is stronger because the starting angle of attack is not null. The grid refinement follows the unsteady behavior of the flow, combining effectively mesh refinement and coarsening, so that grid spacing is reduced in the regions exhibiting Mach variations. 
In addition, mesh refinement near the leading and trailing edge is preserved during the simulation to accurately represent the airfoil geometry.

\begin{figure*}
\centering
\includegraphics[width=0.65\textwidth]{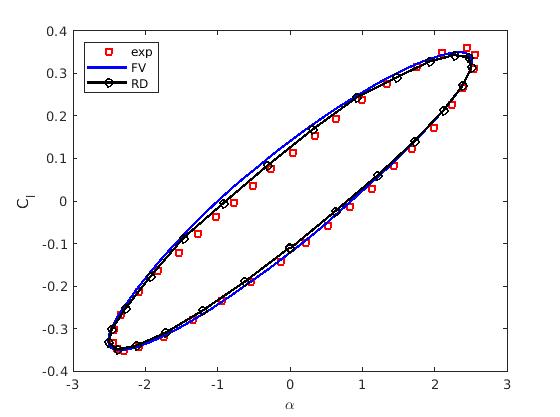} 
\caption{Results of the NACA0012 pitching airfoil: lift coefficient $C_L$ versus angle of attach $\alpha$. A comparison between the experimental data~\cite{Anderson1988}, finite volume and residual distribution scheme is provided.}
\label{f:nacaCLalpha}
\end{figure*}

Finally, Fig.~\ref{f:nacaCLalpha} presents the $C_L-\alpha$ curve, that is the variation of the lift coefficient $C_L$ versus the angle of attack $\alpha$. The results of the proposed ALE RD schemes are compared to the finite volume results~\cite{Re2017}, and the experimental results given in~\cite{Anderson1988}, achieving a fair match.

\section{Conclusions}\label{s:concl}
In this work, we proposed an ALE RD scheme to solve the unsteady Euler equations over adaptive unstructured grids while enforcing automatically the so-called GCL constraint.
First, we derived the ALE formulation of the governing equations for a dynamic fixed-connectivity grid following an analogy with the previous work on finite volume schemes.
Then, we formulated a local expansion-collapse procedure to describe node insertion, node deletion, and edge swap in terms of fictitious continuous deformations. This procedure aims to define ALE residual and GCL-compliant interface velocities that takes into account the local grid modifications, but that can be included straightforwardly into the ALE formulation of the governing equations.
Any explicit interpolation of the solution from the original to the adapted grid is thus avoided, and the numerical properties of the underlying fixed-connectivity ALE RD scheme are preserved. 

The presented work represents a first step toward the more ambitious goal to develop a high-order, conservative and interpolation-free scheme for the simulation of unsteady flow problems involving complex geometries and/or around moving bodies, for which mesh adaptation is an indispensable tool. The work and the 2D results presented here show the goodness and feasibility of the idea. Clearly, further developments are needed to reach the long-term target, such as the increase of the order of accuracy and a possible inclusion of $p$-adaptation techniques (locally adaptive order of the spatial scheme), the extension toward the 3D, which requires also a more efficient software implementation, and the treatment of topological modifications of the domain, as, for instance, the separation or coalition of solid bodies.

\section*{Acknowledgement}
This work was initiated while the authors were at the Institute of Mathematics at University of Z\"{u}rich (UZH). The authors would like to express their gratitude to prof. R\'{e}mi Abgrall for his valuable and constructive suggestions during the planning and development of this research.
S. Colombo would like to thank him also for the financial support during his stay in Z\"{u}rich.
B. Re  gratefully acknowledges the financial support received under the grant \textit{Forschungskredit} of the University of Z\"{u}rich, grant no. [FK-20-121].

\bibliography{biblio}

\end{document}